# ANFIS and metaheuristics for green supply chain with inspection and rework

Nidhi Sharma[1], Madhu Jain[2], Dinesh Sharma[3]
Indian Institute of Technology Roorkee[1,2], Roorkee 247667 (India)
University of Maryland Eastern Shore, USA[3]
n_sharma@ma.iitr.ac.in[1] , madhu.jain@ma.iitr.ac.in[2] , dksharma@umes.edu[3]

**Abstract**
The focus of present article is to investigate a supply chain inventory model of deteriorated items along with inspection and stock dependent demand using green technology to reduce carbon emissions. Products that are decaying have a high sensitivity to the environment in terms of temperature, carbon emission, humidity, waste disposal, etc. This study develops a profit maximization model in the presence of deterioration, preservation, imperfect production, inspection error, rework, stock and price-dependent demand. Three carbon emission strategies are proposed to reduce the expenses in different carbon emissions scenarios. The suggested approach may be used to determine the optimal production period, preservation investment, and level of green investment. The solution of the proposed non-linear constraint optimization is provided by using a penalty method in metaheuristic approaches. In order to conduct a sensitivity analysis for the essential model parameters, a numerical example is presented. The results produced by DE and PSO are compared with the results obtained by Adaptive Neuro-Fuzzy Inference System (ANFIS) technique.

## 1. Introduction

The designing, planning, controlling, and overseeing supply chain (SC) operations used to generate net value, utilizing global logistics, synchronizing supply with demand, building a competitive infrastructure, and evaluating global performance is referred as supply chain management (SCM). The climatic and geopolitical uncertainty, energy shortages, the rising cost of living, unreliability, and lack of transparency etc., are the major problems faced by global SC businesses. Now-a-days, the companies are increasing their green initiatives to offset the effects of conventional operating procedures and carbon emissions as users so as to voluntarily associate themselves with brands that are more accountable for the sustainable development. As mentioned in Coady et al. (2015), an important new estimate by the International Monetary Fund (IMF) has been analysed that the fossil fuel sector benefits from worldwide subsidies of $5.3 trillion a year, or $10 million a minute every day. A green SC can limit the emission of greenhouse gases like Methane ($CH_4$) and Carbon Dioxide ($CO_2$) by directly lowering down the use of fossil fuels while minimizing their carbon footprint. One-third of all anthropogenic greenhouse gas (GHG) emissions come from the food systems. According to estimates from the Food and Agriculture Organization (FAO) of the United Nations, growing emissions by agricultural production reached 10.7 Gt $CO_2$e per year (FAO, 2019). According to Menegat et al. (2022), the SC for synthetic N fertilizer emitted $129.1 \pm 171.1$ (mean $\pm$ s.d) mega tonnes of $CO_2$ in 2018. It is possible to significantly cut emissions by reducing the overall production and usage of synthetic fertilizers. We can promote organic fertilizers and food production systems by reducing $CO_2$ emissions.

The deterioration is described as a loss in usefulness from the original condition due to decay, evaporation, loss of marginal value of the product, such as prescription drugs, dairy products, food items, fruits, vegetables, and electronics. The food business is concerned about food wastage during storage since it is a significant environmental issue which significantly affects SC profitability. Delivering faulty products to clients can lead a decline in reputation, or immensely high expenses. In order to preserve the quality and quantity of products, inventory items, components, and goods should be inspected regularly. Duffuaa and Khan (2005) presented an inspection plan for the quality control of the items that addressed various inspection errors (IE), performance measurements, and the financial impact of the various IE. It is observed that deterioration affects the shelf life of products, business profits, the environment, consumer satisfaction, and especially the cost of items; therefore, we must minimize this loss by using preservation technology (Das et al., 2020).

Soft computing techniques can be used for the grouping, membership, and categorization of many variables that occur in real-world situations. It differs from conventional computing methods because it accepts ambiguity, imperfection, and incomplete truth. The common soft computing approaches are evolutionary computing, artificial neural networks (ANN), fuzzy logic. Together, they can be utilized to address issues that are too complex or inherently noisy to be resolved by employing conventional mathematical techniques. The highly non-linear optimization problems can be easily solved using particle swarm optimization (PSO) and differential evolution (DE). Kennedy and Eberhart (1995) developed the PSO method by considering the behavior of particles in a swarm. The primary inspiration for PSO was the swarming behavior of birds. PSO uses progressively improved solutions based on particle mobility and their interactions to optimize complex problems. Unlike conventional optimization techniques, PSO does not need the objective function to be differentiable. The DE metaheuristic was first conceptualized by Storn and Price (1995) that can tackle the non-linear, multimodal, and non-differentiable optimization issues. Jain and Singh (2022) used genetic algorithm (GA) and DE algorithms with various mutation strategies to resolve the optimization problem of inventory control by considering the price-sensitive demand, inspection sensitive deterioration and partial advance payment. Some researchers have considered the optimization issues of inventory system by considering the price-sensitive demand and inspection-sensitive degradation, and partial advance payment. Jain et al. (2022) outlined some crucial supply chain domains where artificial intelligence (AI) can be used. Akhtar et al. (2023) evaluated the total profit using a unique hybrid approach based on DE and social group optimization methodologies. Lagaros et al. (2023) studied the use of specialized metaheuristic approaches that simplify the handling of performance and bound constraints of optimization problems. To maximize profit for a degraded item, Khedlekar & Kumar (2024) used PSO and worked to determine when, how much, how long, and how much to convert from one type to another form. In multi-echelon distribution systems for perishable commodities under a restricted number of evaluations, Liu & Nishi (2024) suggested an evolutionary optimization technique.

Soft computing techniques like AI can be used for the revolutionizing logistics and SC management by anomaly detection, end-to-end visibility, intelligent decision-making, scheduling maintenance, enhanced customer service, streamlined inventory management, reducing manual work, demand forecasting, fraud prevention, delivery prediction, cost optimization, and dynamic-real time route optimization. A knowledge-based expert system called PILOT was created by Anagun (1997) to analyse an inventory model. In order to explain the uncertainty in demand and some expenses, we are using adaptive neuro-fuzzy inference systems (ANFIS) in this article. As a learning and predictive model for performance evaluation of time and price, the ANFIS model can be applied on training data in carbon taxation policy. Fuzzy parameters and artificial neural networks (ANN) via ANFIS bring the suggested model closer to a realistic approach and offer insightful information for the future digital SC architecture. Güneri et al. (2011) worked on an analytical method for supplier selection decision-making that was based on input selection using the ANFIS model. Sremac et al. (2018) concluded that the hybrid ANFIS is an important concept for strengthening the calculating economic order quantity (EOQ) of a logistic system. In order to effectively distribute products, Okwu et al. (2023) emphasized the necessity for a hybrid intelligent approach and implemented ANFIS in a double source multi-destination system.

The present study on the centralized two-tier green SC between a single producer and retailer is expected to explore the following key research questions:

- Regarding the combined profit and sustainability of the centralized two-tier SC, what is the effect of integrating degradation, preservation, and rework?
- How do type-I and type-II human-based inspection errors affect the overall profitability and efficiency of the SC?
- What is the relationship between the combined profit and environmental effect of the SC and the various carbon emission regulation schemes (carbon taxes, cap & trade, and restricted emission)?
- How does the retailer's goodwill impairment investment affect the total profitability of the SC, and what function does it play in preserving consumer satisfaction?
- What is the impact of incorporating modern techniques like ANFIS on improving the SC inventory model's practical applicability and optimization?

Rework, degradation, preservation, inspection errors, price and stock-dependent demand, carbon emissions, and green investment are all incorporated into a SC inventory model that integrates a single producer and store in order to investigate these research problems. Type-I and Type-II inspection errors that happen during the manufacturer's inventory cycle were also included in this model. We considered the retailer's goodwill impairment investment in order to preserve customer satisfaction, and fully backlogged shortages are permitted at the retailer's level.

It is worth write to review the literature that addresses a few critical elements related to the SCM model developed in the present article. Firstly, we emphasize the notable contributions to green supply chains and present the literature about sustainability and green investment (GI). The articles that make it possible to gain insightful knowledge on preventing deterioration using preservation have been mentioned. We also review the contributions which involve the inspection to enhance the quality of SC. Furthermore, literature review is conducted on the price and stock-dependent variable demand. In the end, we mention some insightful works in the area of SC along with some research gap.

## *1.1.* Inventory models and green investment in SC

The GI play vital role in the reduction of emissions, operating costs, customer loyalty, boosting returns on used goods, and enhancing SC environmental performance. Hua et al. (2011) worked on carbon footprints management in SC management. Zhang & Xu (2013), Qin et al. (2015), Xia et al. (2018), Taleizadeh et al. (2021), Astanti et al. (2022), and many other researchers have worked on carbon taxation, cap and trade carbon emission (CE) and GI policies in SC inventory models. Huang et al. (2020) looked at a two-echelon SC with three CE policies related to restricted total CE, carbon taxes, and cap and trade. In a green production inventory model, Ruidas et al. (2022) evaluated the effects of combined investments in greening innovation and emission reduction technologies by considering a selling price dependent demand. Mala et al. (2022) and Marchi & Zanoni (2023) worked on carbon taxes, cap and trade, and emission limit regulations while considering sustainability, green investment, and logistics. Abbasi & Ahmadi Choukolaei (2023) focused on different approaches and models to examine the impact of carbon laws on green SC network design after reviewing the literature from 2010 to 2023. Jauhari et al. (2023) suggested a single vendor and single buyer model to calculate the frequency of shipments, review duration, safety factor, and GI in order to minimize the combined total cost and emissions in the SC. Pervin (2024) suggested a sustainable inventory approach for items to reduce carbon emissions.

## *1.2.* SC inventory models and deterioration

Items that deteriorate over time and whose actual volume may change, such as gasoline, lubricants, radioactive materials, glues, and chemicals, cannot be stored for a long duration. Wee et al. (2005) used deteriorating rate as Weibull distribution with two parameters for an inventory model and two warehouse inventory models. In order to determine the optimal joint total cost in SC that included the supplier, manufacturer, and consumer, Rau et al. (2003) created a multi-echelon inventory model for a deteriorating item. In order to find a cost-effective approach with an integer number of deliveries and the optimum lot size for three different models, Sarkar (2013) employed probabilistic deterioration. Sarkar et al. (2016) worked on an inventory model by considering the deterioration rate inversely proportional to the reliability. Tiwari et al. (2018) created an inventory model with permissible shortages for deteriorating products under a two-level partial trade credit policy. Under the context of SC integration, Shen et al. (2019) looked into a production inventory model for deteriorating goods under a carbon tax policy and joint preservation technology investment. Using consignment stock and vendor-managed inventory, Hemmati et al. (2023) created a multi-echelon model where demand of a deteriorated item is based on both stock and price. In order to characterize the structure of the optimal inventory policy, Ghasemzadeh & Pamucar (2023) addressed the management of deteriorating commodities of dairy industry in a three-echelon SC network with the help of finite-horizon semi-markov process and metaheuristics for optimization purpose. In order to decrease the environmental effect and increase SC profitability and efficiency, San-José et al. (2024) worked on a sustainable inventory model of deteriorating products that followed a power demand pattern.

*1.3. SC inventory models and inspection errors*
The seller inspects the items to ensure their quality prior to the delivery. The following models explain the different kind of errors caused by human negligence in inspection. Taheri-Tolgari et al. (2012) assumed that the first-stage inspector of product quality control might make some IE during the separation of faulty and perfect items, and following the rework process when there were no inspection faults. Pal & Mahapatra (2017) proposed a three-layer SC production inventory model under the assumption that the manufacturer produces both perfect and imperfect quality product. There may be possibility of IE by labelling non-defective items as defective or defectives as non-defectives. Cheikhrouhou et al. (2018) discussed an inventory model with two possible cases for the product returns. In the first scenario, defective lots were immediately removed from the system and sent back to the supplier from the retailer. In the other one, the retailer sent defective lots back to supplier when receiving the subsequent lot from the supplier. Manna et al. (2020) developed an imperfect production inventory model with selling price discount and warranty-dependent demand under the consideration of IE and time-related cost. Sepehri and Gholamian (2023) explained that shortages can affect a sustainable inventory system of imperfect quality goods when quality improvement and inspection processes are considered concurrently. Chandramohan et al. (2023) considered goods of varying quality, non-instantaneous deterioration, learning effects, and multiple credit periods, together with a carbon tax and retailer-end inspection to preserve consumer goodwill. Wang et al. (2024) investigated the effects of warranty charges and inspection errors on decisions about pricing and quality under several scenarios.

*1.4. SC inventory models with price and stock-dependent demand*
To reduce waste and prevent shortages caused by stock, pricing, time, and a variety of other factors, every business owner tries to forecast customer demand. Inventory management may be influenced by the price and stock-dependent demand. As prices rise, the demand may also change depending on how an item is useful. Some researchers (Urban 1992; Kevin Weng 1995; Bhunia & Maiti 1997; Mondal et al. 2003; Wu et al. 2006; Kilian & Park 2009; Chen et al. 2010 and Chang 2013) used price-dependent demand together with a number of other considerations that may have an impact on the business's profitability. Ghoreishi et al. (2015) worked on an economic ordering policy for non-instantaneously degrading products with selling price and inflation-induced demand under the influence of payment delays, customer returns, and shortages. Mishra et al. (2017) designed an EOQ inventory model that considered two distinct demand rate functions dependent on stock and selling price. Das et al. (2020) created an inventory model for non-instantaneous depreciating goods, where the selling price influenced the demand. To address highly non-linear optimization issue, quantum-behaved PSO (QPSO) variations was utilized. Abdul Halim et al. (2021) suggested a production inventory model for perishable goods with non-linear pricing, linear stock-dependent market demand, and a strategy for overtime production. In order to incorporate the ecological initiatives, Shah et al. (2023) developed a perishable inventory model by considering stock, price and greening level dependent demand. With various warehouses and price- and stock-dependent demand, Rodríguez-Mazo & Ruiz-Benítez (2024) developed a deterministic replenishment model.

## *1.5.* **Research gap**

The suggested study on green SCM is done by identifying the gaps after surveying the relevant literature. Inspection not only affects the product's quality, but also contributes to the future sales. Product inspection is essential for manufacturing sectors to produce perfect goods and for the retailers to maintain goodwill. The manufacturers and retailers generally look upon a balance between the quantity and price of a product with growing sales. To the author's best knowledge, no research has been conducted on a green SC inventory system that is connected to rework, deterioration, preservation all together. The carbon policies, and IEs in which demand is stock and price-dependent are to be involved. To highlight the noble features of proposed model, a comparison of different characteristics used in the literature and current research work is presented in Table 1.

    The subsequent sections of the article are organized as follows. Section 2 provides a brief discussion of the model formulation along with requisite notations and assumptions. Section 3 presents the mathematical analysis used to frame the total profit function for the manufacturer and retailer. Section 4 is devoted to frame the profit function in three different CE scenarios. The working of metaheuristic optimization, ANFIS results and concavity are explained in Section 5 by using numerical simulations. Section 6 focuses on sensitivity testing of the input parameters. Managerial implications are covered in Section 7. The last Section 8, concludes the entire research output.

**Table 1.** A comparative analysis of the proposed and previous studies

| S. No | Author(s) | Demand function | IE | Deterioration | Goodwill Impairment | Environmental impact/ CE | Rework | Meta-Heuristics | ANFIS |
|---|---|---|---|---|---|---|---|---|---|
| 1. | Duffuaa & Khan, (2005) | Constant | Yes | No | No | No | Yes | No | No |
| 2. | Chaharsooghi et al. (2008) | Constant | No | No | No | No | No | GA & RL | No |
| 3. | Taheri-Tolgari et al. (2012) | Constant | Yes | No | No | No | Yes | No | No |
| 4. | Ghoreishi et al. (2015) | Price dependent | No | Yes | No | No | No | No | No |
| 5. | Sremac et al. (2018) | Constant | No | No | No | No | No | No | No |
| 6. | Shen et al. (2019) | Constant | No | Yes | No | Yes | No | No | No |
| 7. | Chandra Das et al. (2020) | Price dependent | No | Yes | No | No | No | PSO, QPSO, AQPSO, GQPSO, WQPSO | No |
| 8. | Hemmati et al. (2021) | Price & stock dependent | No | Yes | No | No | No | Nelder–Mead Algorithm | No |
| 9. | Ruidas et al. (2022) | Price & greenness level dependent | No | No | No | Yes | Yes | QPSO | No |
| 10. | Jauhari et al. (2023) | Stochastic | No | No | No | Yes | No | No | No |
| 11. | San-José et al. (2024) | Power | No | Yes | No | Carbon tax | No | No | No |
| 12. | **This paper** | **Price & stock dependent** | **Yes** | **Yes** | **Yes** | **Carbon taxation, Cap & trade, and Limited emission policies** | **Yes** | **DE & PSO** | **Yes** |

## 2. Model formulation

In this study, an integrated single manufacturer and retailer inventory model with fully backlogged shortages, isformulated. The present article incorporates many realistic features such as (i) inspection space, (ii) deterioration, (iii) preservation technology, (iv) stock and price dependent demand, (v) goodwill impairment, (vi) carbon emission, (vii) green technology, (viii) carbon tax, (ix) cap & trade and (x) CE under a limited capacity constraint.

### 2.1. Notations and assumptions

Following are some pertinent notations that have been used to formalize the proposed model:

**Parameters:**

| | | |
|---|---|---|
| $T_1$ | : | Rework time for the defectives produced at manufacturer's level (years) |
| $T_{11}$ | : | Retailer's time in which shortages are fulfilled (years) |
| $T_2$ | : | Time needed for the retailer to build inventory (years) |
| $T_3$ | : | Complete cycle time of the retailer (years) |
| $P, (P_e)$ | : | Production (effective production) rate of the manufacturer (units/year) |
| $P_r$ | : | Rework rate at the manufacturer's level (units/year) |
| $P_{de}$ | : | Effective production rate for defectives (units/year) |
| $\beta_1, (\beta_2)$ | : | Probability of type- I (type- II) error in screening process at the manufacturer's level |
| $f_d$ | : | Proportion of generated defective at the manufacturer level |
| $\theta_1, (\theta_2)$ | : | Deterioration rate for the manufacturer (retailer) |
| $Q_m, (Q_r)$ | : | Economic order quantities for the manufacturer (retailer) (units) |
| $D_r$ | : | Demand rate for the retailer (units/year) |
| $C_p, (C_r)$ | : | Production (rework) cost for the manufacturer ($/unit) |
| $C_g$ | : | Penalty costs for uninspected defective items ($/unit) |
| $C_{op}, (C_{or})$ | : | Setup costs for production (rework) process ($/setup) |
| $i_c$ | : | Inspection cost for the manufacturer ($/unit) |
| $h_p, (h_d)$ | : | Holding cost for perfect (defective) items for the manufacturer ($/unit/year) |
| $h_r, (C_s)$ | : | Holding (shortages) cost for the retailer ($/unit/year) |
| $d_{cp}, (d_{cd})$ | : | Deterioration cost for perfect (defective) items for the manufacturer ($/unit) |
| $a, (b)$ | : | Constant (price dependent) parameter for customer's demand |
| $W_m$ | : | Selling price for the manufacturer ($/unit) |
| $C_{dr}, O_r$ | : | Deterioration cost ($/unit) and ordering cost ($/order) for the retailer |
| $E_p, (E_t)$ | : | Emission cost (EC) due to production (transportation) of the goods ($/unit) |
| $E_{h1}, (E_{h2})$ | : | EC due to holding the perfect (defective) goods ($/unit/year) |
| $E_{\theta_1}, (E_{\theta_2})$ | : | EC due to deterioration from the perfect (defective) goods ($/unit) |
| $E_{h3}, (E_{\theta_3})$ | : | EC due to holding (deterioration) of the retailer ($/unit/year) |
| $d_1$ | : | Travelled distance to supply the retailer's order at manufacturer's level (Kilometres) |
| $C_{Tax}$ | : | Carbon tax per unit CE ($/tonne) |
| $C_{C\&T}$ | : | Expected cost at which trade of carbon will be done ($/tonne) |
| $U_1, (U_2)$ | : | Maximum limits of CE in cap & trade (limited emission policies) (Tonne) |

| $s$ | : | Shortage quantity at the retailer's level (Units) |
| $\kappa_1, (\kappa_2)$ | : | Exponents of carbon reduction functions for carbon taxation (limited emissions) ($\kappa_1, \kappa_2 > 0$) |

**Decision variables:**

| $T_0$ | : | Effective production time (years) |
| $\xi_1, \xi_2$ | : | Preservation cost per unit item per unit time for the manufacturer and retailer ($/year) |
| $G$ | : | Green investment parameter ($/year) |
| $W_r$ | : | Selling price for the retailer ($/unit) |

**Abbreviation:**

| SC | : | Supply chain |
| SCM | : | Supply chain management |
| GSCM | : | Green supply chain management |
| CE | : | Carbon emission |
| GI | : | Green investment |
| GHG | : | Greenhouse gas |
| PSO | : | Particle swarm optimization |
| DE | : | Differential evolution |
| EC | : | Emission cost |

### 2.2. Assumptions

To study a centralized two-tier green SC between a single manufacturer and retailer, the following assumptions are made:

a) At manufacturer's level, there will be faulty products since it is assumed that the production process is imperfect. Rework procedure is used on defectives to reduce waste.

b) The quality checker can make one of two sorts of mistakes during the examination. The likelihood that a non-faulty item will be discarded as defective, is called human-based type-I IE ($\beta_1$). When a faulty item is declared as non-defective, this is known as human-based type-II IE ($\beta_2$). These errors are occurring at the manufacturer's inventory cycle.

c) The demand function ($D_c$) for the customer is dependent on the selling price ($W_r$) and instantaneous stock level ($I(t)$) which is given as follows:

$$D_c(W_r, I(t)) = \begin{cases} f(W_r) + \eta I(t), & I(t) > 0 \\ f(W_r), & I(t) \leq 0 \end{cases} \quad (1)$$

where $\eta$ is stock dependent consumption parameter lying in $(0,1]$ and $f(W_r) = a - bW_r$.

d) At the retailer's level, some defectives may arrive due to IEs at the manufacturer's level.

e) The retailer will invest $f_r\%$ of the total profit in goodwill impairment to maintain the customer satisfaction and goodwill in the market.

f) Completely backlogged shortages are allowed at the retailer's level.

g) In the SC system, preservation technology is utilized to reduce the impact of degradation at inventory level of both chain members. The preservation technology functions for the manufacturer and retailer are $m_1(\xi_1) = 1 - e^{-v_1\xi_1}$ and $m_2(\xi_2) = 1 - e^{-v_2\xi_2}$ respectively, where $v_1, v_2 > 0$.

h) CE is considered during the processes of manufacturing, storage, decay, and transportation.
i) The three strategies viz., (i) emission tax reduction strategy, (ii) cap & trade strategy and (iii) restricted CE strategy have been employed to assess the environmental impacts.
j) The manufacturer is paying '$\omega$' fraction of GI and the remaining part is paid by the retailer in whole SC. When using green technology, CE for each player are reduced by $\rho_m = \omega G l_1 - l_2(\omega G)^{\kappa_1}$ and $\rho_r = (1-\omega)G l_1 - l_2((1-\omega)G)^{\kappa_1}$, where $\rho_m$ and $\rho_r$ stand for the manufacturer and the retailer and $l_1, l_2$ are the carbon reduction efficiency factor and offsetting carbon reduction factors, respectively.

## 3. Profit analysis

The mathematical formulations for the profit function for the manufacturer and retailer are given as follows:

### 3.1. The manufacturer's profit function

In centralized SC, a single product is produced at a finite production rate $P$ by the producer after receiving the retailer's demand. The manufacturer produces both perfect and flawed goods at rates of $(1 - f_d)P$ and $f_d P$ respectively, throughout each production run. The screening procedure is carried out in accordance with the manufacturing cycle. Two types of errors $\beta_1$ and $\beta_2$ are arising in screening process. We can define the loss of perfect units and addition of imperfect units in lot due to type-I and type-II errors as $(1 - f_d)P\beta_1$ and $f_d P \beta_2$, respectively. Thus, the effective production rate up to time $T_0$ is given by

$$P_e = (1 - f_d)P + f_d P \beta_2 - (1 - f_d)P\beta_1 \tag{2}$$

The manufacturer repairs the defectives at rate $P_r$ from time $T_0$ to $T_1$. In order to decrease degradation during the whole manufacturing cycle, we consider the preservation investment following the rate $m_1(\xi_1) = 1 - e^{-v_1\xi_1}$ where $v_1 > 0$. The inventory level of the manufacturer is shown in Fig. 1.

The rate of deterioration for the manufacturer is given as $\theta_m = \theta_1(1 - m_1(\xi_1))$.

In the manufacturer cycle, over the period $[0, T_2]$, the differential equations given below determine the inventory levels:

$$\frac{dI}{dt} = \begin{cases} (1 - f_d)P + f_d P \beta_2 - (1 - f_d)P\beta_1 - \theta_m I & ; 0 \leq t \leq T_0 \\ P_r - \theta_m I & ; T_0 \leq t \leq T_1 \\ -D_r - \theta_m I & ; T_1 \leq t \leq T_2 \end{cases} \tag{3}$$

with initial and boundary conditions $I(0) = 0; I(T_1) = Q_m$ and $I(T_2) = 0$.

The solution of eq. (3) w.r.t. the initial and boundary conditions is given as follows:

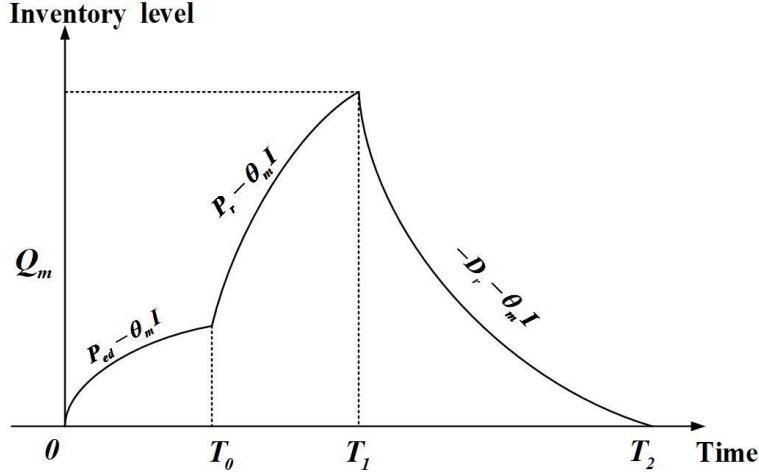

**Fig 1**: Inventory flow of the manufacturer

$$I(t) = \begin{cases} \frac{P_e}{\theta_m}(1 - e^{-\theta_m t}) & ; 0 \leq t \leq T_0 \\ \frac{P_r}{\theta_m} + \left(Q_m - \frac{P_r}{\theta_m}\right)e^{\theta_m(T_1-t)} & ; T_0 \leq t \leq T_1 \\ \frac{D_1}{\theta_m}\left(e^{\theta_m(T_1-t)} - 1\right) & ; T_1 \leq t \leq T_2 \end{cases} \quad (4)$$

By continuity conditions at $t = T_0$ and $t = T_1$, we get $T_1$ and $Q_m$ as given below:

$$T_1 = T_0 + \frac{1}{\theta_m}\left(\frac{P_e(1-e^{-\theta_m T_0})-P_r}{Q_m \theta_m - P_r}\right) \quad (5)$$

$$Q_m = \frac{D_r}{\theta_m}\left(e^{\theta_m(T_2-T_1)} - 1\right) \quad (6)$$

Inventory level for defective items by using effective defective production rate ($P_{de}$), can be obtained by using following differential equations where $P_{de} = (1 - f_d)P\beta_1 + f_d P - f_d P\beta_2$.

$$\frac{dI_d}{dt} = \begin{cases} P_{de} - \theta_m I_d & ; 0 \leq t \leq T_0 \\ P_r - \theta_m I_d & ; T_0 \leq t \leq T_1 \end{cases} \quad (7)$$

with initial and boundary conditions as $I_d(0) = 0$; $I_d(T_1) = 0$.

The solution of (7) by using initial and boundary conditions is given below:

$$I_d(t) = \begin{cases} \frac{P_{de}}{\theta_m}(1 - e^{-\theta_m t}) & ; 0 \leq t \leq T_0 \\ \frac{P_r}{\theta_m}\left(e^{\theta_m(T_1-t)} - 1\right) & ; T_0 \leq t \leq T_1 \end{cases} \quad (8)$$

By continuity conditions at $t = T_0$, we get $T_1$ by using eq. (8) as

$$T_1 = T_0 + \frac{1}{\theta_m} \log\left(1 + \frac{P_{de}(1-e^{-\theta_m T_0})}{P_r}\right) \quad (9)$$

By using (5) and (9), we get

$$Q_m = \frac{P_r}{\theta_m}\left(1 + \frac{P_e(1-e^{-\theta_m T_0})-P_r}{P_{de}(1-e^{-\theta_m T_0})+P_r}\right) \quad (10)$$

By using eq. (6), we get

$$T_2 = T_1 + \frac{1}{\theta_m} \log\left(1 + \frac{Q_m \theta_m}{D_r}\right) \tag{11}$$

Cost components associated with the manufacturer's SC activities are established as given below:

- Sales revenue:
$$SR_m = W_m D_r (T_2 - T_1) \tag{12.a}$$

- Production cost:
$$PC_m = C_p P T_0 \tag{12.b}$$

- Setup cost:
$$StC_m = C_{op} + C_{or} \tag{12.c}$$

- Penalty cost:
$$PeC_m = C_g f_d P \beta_2 T_0 \tag{12.d}$$

- Rework cost:
$$RC_m = C_r P_r (T_1 - T_0) \tag{12.e}$$

- Preservation cost:
$$PreC_m = \xi_1 T_2 \tag{12.f}$$

- Screening cost:
$$ScC_m = i_c P T_0 \tag{12.g}$$

- Holding cost for perfect items:
$$HC_{m1} = h_p \int_0^{T_2} I(t)dt = h_p \left[\int_0^{T_0} I(t)dt + \int_{T_0}^{T_1} I(t)dt + \int_{T_1}^{T_2} I(t)dt\right]$$

$$= h_p \left\{\frac{P_e}{\theta_m}\left(T_0 + \frac{1}{\theta_m}\left(e^{-\theta_m T_0} - 1\right)\right) + \frac{P_r}{\theta_m}(T_1 - T_0) + \left(Q_m - \frac{P_r}{\theta_m}\right)\frac{(e^{\theta_m(T_1 - T_0)} - 1)}{\theta_m} + \frac{D_r}{\theta_m}(T_1 - T_2) + \frac{D_r}{\theta_m^2}\left(e^{\theta_m(T_2 - T_1)} - 1\right)\right\} \tag{12.h}$$

- Holding cost for defective items:
$$HC_{m2} = h_d \int_0^{T_1} I_d(t)dt = h_d \left[\int_0^{T_0} I_d(t)dt + \int_{T_0}^{T_1} I_d(t)dt\right]$$

$$= h_d \left\{\frac{P_{de}}{P_e}\left(\frac{P_e}{\theta_m}\left(T_0 + \frac{1}{\theta_m}\left(e^{-\theta_m T_0} - 1\right)\right)\right) + \frac{P_r}{\theta_m}(T_0 - T_1) + \frac{P_r}{\theta_m^2}\left(e^{\theta_m(T_1 - T_0)} - 1\right)\right\} \tag{12.i}$$

- Deterioration cost for perfect items:
$$DC_{m1} = d_{cp} \int_0^{T_2} \theta_1 I(t)dt = d_{cp} \theta_1 \left(\frac{HC_{m1}}{h_p}\right) \tag{12.j}$$

- Deterioration cost for defective items:
$$DC_{m2} = d_{cd} \int_0^{T_1} \theta_1 I_d(t)dt = d_{cd} \theta_1 \left(\frac{HC_{m2}}{h_d}\right) \tag{12.k}$$

The manufacturer's total profit without accounting for expenditures associated with CE is

$$\varphi_m = \frac{1}{T_2}[SR_m - (PC_m + StC_m + PeC_m + RC_m + PreC_m + ScC_m + HC_{m1} + HC_{m2} +$$

$$DC_{m1} + DC_{m2})] \tag{13}$$

**CE related cost components at the manufacturer level**

When considering the CE from the processes of degradation, production, transportation and storage with GI, the quantity of CE for manufacturer is obtained as

- EC due to production:
$$e_{m1} = Q_m E_p \tag{14.a}$$
- EC due to holding perfect items:
$$e_{m2} = E_{h1}\left(\frac{HC_{m1}}{h_p}\right) \tag{14.b}$$
- EC due to holding defective items:
$$e_{m3} = E_{h2}\left(\frac{HC_{m2}}{h_d}\right) \tag{14.c}$$
- EC due to deteriorated perfect items:
$$e_{m4} = E_{d1}\left(\frac{DC_{m1}}{d_{cp}}\right) \tag{14.d}$$
- EC due to deteriorated defective items:
$$e_{m5} = E_{d2}\left(\frac{DC_{m2}}{d_{cd}}\right) \tag{14.e}$$
- EC due to transportation:
$$e_{m6} = d_1 E_t D_r (T_2 - T_1) \tag{14.f}$$

The total CE cost is

$$CarC_m = Q_m E_p + E_{h1}\left\{\frac{P_e}{\theta_m}\left(T_0 + \frac{1}{\theta_m}\left(e^{-\theta_m T_0} - 1\right)\right) + \frac{P_r}{\theta_m}(T_1 - T_0) + \left(Q_m - \frac{P_r}{\theta_m}\right)\frac{(e^{\theta_m(T_1 - T_0)} - 1)}{\theta_m} + \frac{D_r}{\theta_m}(T_1 - T_2) + \frac{D_r}{\theta_m^2}\left(e^{\theta_m(T_2 - T_1)} - 1\right)\right\} + E_{h2}\left\{\frac{P_{de}}{P_e}\left(\frac{P_e}{\theta_m}\left(T_0 + \frac{1}{\theta_m}\left(e^{-\theta_m T_0} - 1\right)\right)\right) + \frac{P_r}{\theta_m}(T_0 - T_1) + \frac{P_r}{\theta_m^2}\left(e^{\theta_m(T_1 - T_0)} - 1\right)\right\} + E_{d1}\theta_1\left\{\frac{P_e}{\theta_m}\left(T_0 + \frac{1}{\theta_m}\left(e^{-\theta_m T_0} - 1\right)\right) + \frac{P_r}{\theta_m}(T_1 - T_0) + \left(Q_m - \frac{P_r}{\theta_m}\right)\frac{(e^{\theta_m(T_1 - T_0)} - 1)}{\theta_m} + \frac{D_r}{\theta_m}(T_1 - T_2) + \frac{D_r}{\theta_m^2}\left(e^{\theta_m(T_2 - T_1)} - 1\right)\right\} + E_{d2}\theta_1\left\{\frac{P_{de}}{P_e}\left(\frac{P_e}{\theta_m}\left(T_0 + \frac{1}{\theta_m}\left(e^{-\theta_m T_0} - 1\right)\right)\right) + \frac{P_r}{\theta_m}(T_0 - T_1) + \frac{P_r}{\theta_m^2}\left(e^{\theta_m(T_1 - T_0)} - 1\right)\right\} + d_1 E_t D_r (T_2 - T_1) \tag{15}$$

### 3.2. Formulation of the retailer's profit function

If the retailer has no item to meet the customer's demand by time $T_1$, it must deal with shortages and price-dependent demand after $T_1$ as shown in Fig. 2. We assume that the shortages are fully backlogged. After time $T_1$, the retailer has inventory to fulfil the customer's demand. The retailer is fulfilling the demand of the customers at the rate $D_r - D_c$ up to time $T_{11}$ without deterioration. After time $T_{11}$, the retailer is building inventory of deteriorated items and satisfy the demand of customers up to time $T_3$. To maintain the degradation, the retailer is using preservation investment following the function $m_2(\xi_2) = 1 - e^{-v_2 \xi_2}$, where $v_2 > 0$.

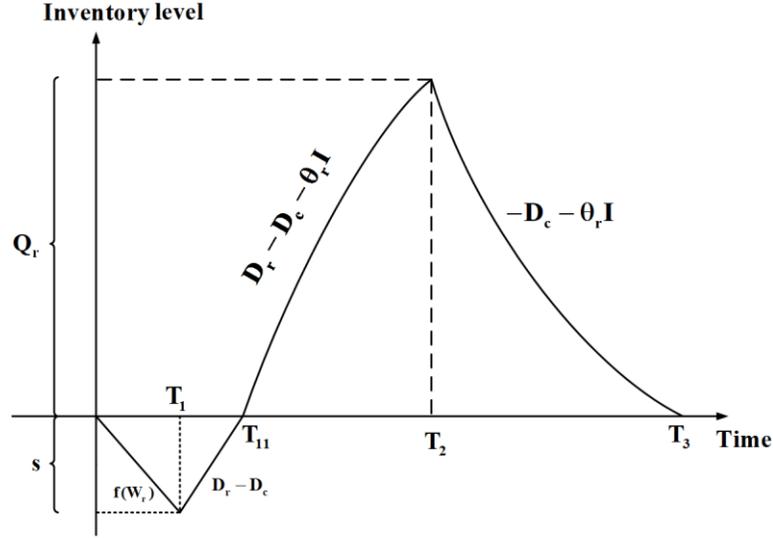

**Fig 2**: Inventory level of retailer

The rate of deterioration for the retailer is given by $\theta_r = \theta_2(1 - m_2(\xi_2))$. In the retailer's cycle, the customer's demand is price and stock dependent to preserve shortages from time $T_{11}$ to $T_3$. The inventory levels by including the effect of shortages, demand and deterioration are governed by the following differential equations.

$$\frac{dI}{dt} = \begin{cases} f(W_r) & ; 0 \leq t \leq T_1 \\ D_r - D_c & ; T_1 \leq t \leq T_{11} \\ D_r - D_c - \theta_r I & ; T_{11} \leq t \leq T_2 \\ -D_c - \theta_r I & ; T_2 \leq t \leq T_3 \end{cases} \quad (16)$$

The initial and boundary conditions for each interval are given by

$I(0) = 0; \ I(T_1) = s; I(T_{11}) = 0; I(T_2) = Q_r$ and $I(T_3) = 0$ \hfill (17)

We denote $B_1 = \eta + \theta_r$ and $B_2 = D_r - f(W_r)$. Solving differential eqs. (16) along with the initial and boundary conditions (17), we get the solutions of equations as:

$s = f(W_r) T_1$ \hfill (18)

$$I(t) = \begin{cases} \frac{B_2}{\eta} + \left(s - \frac{B_2}{\eta}\right) e^{\eta(T_1-t)} & ; T_1 \leq t \leq T_{11} \\ \frac{B_2}{B_1}\left(1 - e^{B_1(T_{11}-t)}\right) & ; T_{11} \leq t \leq T_2 \\ \frac{f(W_r)}{B_1}\left(e^{B_1(T_3-t)} - 1\right) & ; T_2 \leq t \leq T_3 \end{cases} \qquad (19)$$

Using the conditions of $I(T_{11}) = 0$ and $I(T_2) = Q_r$, we get

$$T_{11} = T_1 - \frac{1}{\eta} \log\left\{\frac{B_2}{B_2 - s\eta}\right\} \qquad (20)$$

$$Q_r = \frac{B_2}{B_1}\left(1 - e^{B_1(T_{11}-T_2)}\right) \qquad (21)$$

By continuity at $T_2$, we get

$$T_3 = T_2 + \frac{1}{B_1} \log\left\{1 + \frac{B_1 Q_r}{f(W_r)}\right\} \qquad (22)$$

Sales revenue and cost components for the retailer are given below:

- Sales revenue:

$$SR_r = P_r \int_{T_1}^{T_3} \left(f(W_r) + \eta I(t)\right) dt = P_r f(W_r)(T_3 - T_1) + P_r \eta \int_{T_1}^{T_3} I(t) dt$$

$$= P_r f(W_r)(T_3 - T_1) + P_r \eta \left[\frac{B_2}{\eta}(T_{11} - T_1) + \frac{1}{\eta}\left(s - \frac{B_2}{\eta}\right)\left(1 - e^{\eta(T_1-T_{11})}\right) + \frac{B_2}{B_1}\left(T_2 - T_{11} + \frac{1}{B_1}\left(e^{B_1(T_{11}-T_2)} - 1\right)\right) + \frac{f(W_r)}{B_1}\left(\frac{1}{B_1}\left(e^{B_1(T_3-t)} - 1\right) - (T_3 - T_2)\right)\right] \qquad (23.a)$$

- Holding cost:

$$HC_r = h_r \int_{T_{11}}^{T_3} I(t) dt = h_r \left[\int_{T_{11}}^{T_2} I(t) dt + \int_{T_2}^{T_3} I(t) dt\right]$$

$$= h_r \left[\frac{B_2}{B_1}\left(T_2 - T_{11} + \frac{1}{B_1}\left(e^{B_1(T_{11}-T_2)} - 1\right)\right) + \frac{f(W_r)}{B_1}\left(\frac{1}{B_1}\left(e^{B_1(T_3-t)} - 1\right) - (T_3 - T_2)\right)\right]$$

$$(23.b)$$

- Deterioration cost:

$$DC_r = d_{cr} \int_{T_{11}}^{T_3} \theta_2 I(t) dt = d_{cr} \theta_2 \left(\frac{HC_r}{h_r}\right) \qquad (23.c)$$

- Purchasing cost:

$$PC_r = W_m D_1 (T_2 - T_1) \qquad (23.d)$$

- Ordering cost:

$$OC_r = O_r \qquad (23.e)$$

- Preservation cost:

$$PreC_r = \xi_2 T_2 \qquad (23.f)$$

- Shortage cost:

$$SC_r = \frac{1}{2} s\, T_1\, C_s \qquad (23.g)$$

The retailer's total profit without accounting for expenditures associated with CE and goodwill is

$$\varphi_{r°} = \frac{1}{T_3}[SR_r - (HC_r + DC_r + PC_r + OC_r + PreC_r + SC_r)] \qquad (24)$$

The retailer will pay goodwill charges as a result of a manufacturer inspection error because the retailer is unaware of the flaws in the ordered products (Gu & Lev, 2011). After paying for goodwill impairment (Ye & Yu, 2024), the profit function for the retailer is

$$\varphi_r = (1 - f_r)\varphi_{r°} \qquad (25)$$

**CE costs at the retailer's level**

When considering the CE from the processes of degradation and storage, the amount of CE for the retailer is evaluated using the following cost elements:

- EC cost due to holding perfect items:

$$e_{r1} = E_{hr}\left(\frac{HC_r}{h_r}\right) \qquad (26.a)$$

- EC cost due to deteriorated perfect items:

$$e_{r2} = E_{dr}\left(\frac{DC_r}{d_{cr}}\right) \qquad (26.b)$$

The total CE cost for the retailer is

$$CarC_r = (E_{hr} + E_{dr}\theta_2)\left[\frac{B_2}{B_1}\left(T_2 - T_{11} + \frac{1}{B_1}\left(e^{B_1(T_{11}-T_2)} - 1\right)\right) + \frac{f(W_r)}{B_1}\left(\frac{1}{B_1}\left(e^{B_1(T_3-t)} - 1\right) - (T_3 - T_2)\right)\right] \qquad (27)$$

## 4. Policies for carbon emission and green technology in centralized SC

When governments wish to offer a clear financial incentive to cut emissions across all sectors, they can implement carbon taxes. They are frequently applied as a general economic policy. In industries where emissions are simpler to detect and track, cap & trade is frequently used. It works well for sectors of the economy where reducing emissions with flexibility is advantageous. Direct regulation that places a cap on emissions for certain industries or organizations are considered in limited emission. These regulations may be created to support one another. We have only worked at carbon taxing policy numerically in order to support the model validation and parameters' impact on the total profit of SC.

In the current model, emissions are measured using a reduction function $\rho_m$ and $\rho_r$ for the manufacturer and retailer, respectively. Both partners have invested in environmentally friendly resources. The resulting expenses for a carbon tax, cap & trade and limited CE (Ruidas et al., 2022) by using green investments are listed below.

### (i) Carbon Taxation

Coal, oil, petrol, and natural gas are all carbon-rich fuels that release GHGs when burned. The resulting climatic disturbance results in extreme weather, including heat waves, flooding, and other effects. Carbon taxes can have more immediate positive effects on the environment while addressing climate change by lowering greenhouse emissions.

We consider the carbon tax ($C_{Tax}$) for unit CE in our current model. Both the producer and retailer make investments in green technology to reduce CE and consequently, lowering the cost of paying the $C_{Tax}$. The manufacturing, shipping, storage, degradation, overall cost of paying the $C_{Tax}$, and expense of GI make up the expenses of the current green SC inventory model. The overall profit and expenses for the manufacturer and retailer are as follows:

$$\varphi_m^{Tax} = \frac{1}{T_2}[SR_m - (PC_m + StC_m + PeC_m + RC_m + PreC_m + ScC_m + HC_{m1} + HC_{m2} + DC_{m1} + DC_{m2} + C_{Tax}(e_{m1} + e_{m2} + e_{m3} + e_{m4} + e_{m5} + e_{m6} - \rho_m) + \omega G T_2)] \quad (28)$$

$$\varphi_{r°}^{Tax} = \frac{1}{T_3}[SR_r - (HC_r + DC_r + PC_r + OC_r + PreC_r + SC_r + C_{Tax}(e_{r1} + e_{r2} - \rho_r) + (1-\omega)G(T_3 - T_{11}))] \quad (29)$$

$$\varphi_r^{Tax} = (1 - f_r)\varphi_{r°}^{Tax} \quad (30)$$

The joint profit in centralized SC is given by

$$\varphi_T^{Tax}(T_0, \xi_1, \xi_2, G, W_r) = \varphi_m^{Tax} + \varphi_r^{Tax} \quad (31)$$

and $T_0, \xi_1, \xi_2, G, W_r > 0$.

### (ii) Cap & Trade

The entire quantity of GHG emissions that can be discharged by factories, power plants and other industrial infrastructures is referred to as the cap or limit. Businesses who increase their emissions over the cap pay taxes. Under a cap & trade programme, companies with low CE can exchange extra emission credits they have with other companies that produce more emissions than they are allowed for it.

The combined CE of the manufacturer and retailer are governed by the cap & trade policy under this approach. In order to comply with the constraints of restricted CE, the firm must spend in costs or acquire permits from others if CE exceed the upper limit $U_1$. Let $C_{C\&T}$ be the cost of buying or selling as well as the price of carbon trading. Furthermore, we assume that there is a sufficient supply of CE permits available for the purchase. The overall cost when taking the cap & trade scenario into account includes the costs of retailer's inventory holding, degradation and the amount invested in green technology. The following equations summarize the entire revenue and expenses incurred by the manufacturer and retailer:

$$\varphi_m^{C\&T} = \frac{1}{T_2}[SR_m - (PC_m + StC_m + PeC_m + RC_m + PreC_m + ScC_m + HC_{m1} + HC_{m2} + DC_{m1} + DC_{m2} + C_{C\&T}(e_{m1} + e_{m2} + e_{m3} + e_{m4} + e_{m5} + e_{m6} - \rho_m - U_1) + \omega G T_2)] \quad (32)$$

$$\varphi_{r^\circ}^{C\&T} = \frac{1}{T_3} \left[ SR_r - (HC_r + DC_r + PC_r + OC_r + PreC_r + SC_r + C_{C\&T}(e_{r1} + e_{r2} - \rho_r - U_1) + (1 - \omega)G(T_3 - T_{11})) \right] \tag{33}$$

$$\varphi_r^{C\&T} = (1 - f_r)\varphi_{r^\circ}^{C\&T} \tag{34}$$

The joint profit in centralized SC is given by

$$\varphi_T^{C\&T}(T_0, \xi_1, \xi_2, G, W_r) = \varphi_m^{C\&T} + \varphi_r^{C\&T} \tag{35}$$

and $T_0, \xi_1, \xi_2, G, W_r > 0$.

**(iii) Limited CE**

The producer and retailer must modify their business practices to comply with the restricted CE when considering the limited CE policy into account with upper limit $U_2$. Both parties can invest in green technology to minimize CE. The existing green SC inventory model has expenses for manufacturing, transportation, storage, deterioration, and GI costs. The total amount of CE can be obtained by adding up all the costs associated with CE, the cost and then comparing the result to the upper limit of CE.

The constrained optimization problem for proposed model in the case of limited CE is formulated as follows:

$$Maximize \quad \varphi_T(T_0, \xi_1, \xi_2, G, W_r) - G \tag{36}$$

$$subject\ to\ e_{m1} + e_{m2} + e_{m3} + e_{m4} + e_{m5} + e_{m6} + e_{m7} + e_{m8} - \rho_G \leq U_2 \tag{37}$$

and $T_0, \xi_1, \xi_2, G, W_r \geq 0$.

Here, $\varphi_T = \varphi_m + \varphi_r$. The profit functions for the manufacturer ($\varphi_m$) and retailer ($\varphi_r$) are given in eqs. (13) and (25), respectively.

Also, $\rho_G = Gl_3 - l_4(G)^{\kappa_2}$ where $l_3$, $l_4$ are the efficiency and offsetting carbon reduction parameters, respectively.

Various nonlinear terms are included in the profit functions $\varphi_m$ and $\varphi_r$. It is difficult to use conventional analytical techniques, which frequently rely on concavity assumptions, on nonlinear and non-concave functions. The optimization problems contain a number of decision variables as well as intricate constraints, such as equality and inequality constraints. It can be difficult to find an analytical solution for these optimization problems in all three policies, viz., the carbon tax, cap and trade, and limited emission. In our model, we are concerned with a constrained optimization model which is difficult to solve by using analytic approaches. The formulized optimization issues have been solved using metaheuristic techniques.

## 5. Numerical simulation

The purpose of numerical simulation is to explore the impact of various CE legislations on the overall profit subject to constraints related to production time, green investments, selling prices and preservation costs using a numerical example in the centralized SC. Metaheuristic algorithms, viz., DE and PSO are used to generate the numerical results so as to validate the suggested model. Table 2 contains the values of default parameters for the computational purpose which we have taken by the literature.

**Table 2:** Values of default parameters

| $P$ | $P_r$ | $d_1$ | $W_m$ | $\kappa_1$ | $\kappa_2$ | $f_r$ | $f_d$ | $\beta_1$ | $\beta_2$ |
|---|---|---|---|---|---|---|---|---|---|
| 7500 | 2500 | 25 | 80 | 1.45 | 0.8 | 0.01 | 0.05 | 0.04 | 0.06 |
| $U_1$ | $U_2$ | $C_s$ | $C_p$ | $C_g$ | $C_r$ | $C_{dr}$ | $h_p$ | $h_d$ | $h_r$ |
| 30 | 120 | 2 | 15 | 4 | 5 | 0.05 | 5 | 3 | 2.1 |
| $\theta_1$ | $\theta_2$ | $E_p$ | $E_t$ | $E_{h1}$ | $E_{h2}$ | $E_{h3}$ | $E_{\theta_1}$ | $E_{\theta_2}$ | $E_{\theta_3}$ |
| 0.15 | 0.1 | 0.15 | 0.11 | 0.12 | 0.1 | 0.14 | 0.13 | 0.15 | 0.12 |
| $D_r$ | $\eta$ | $O_2$ | $i_c$ | $(a, b)$ | $(l_1, l_2)$ | $(l_3, l_4)$ | $\omega$ | $d_{cp}$ | $d_{cd}$ |
| 400 | 1.6 | 130 | 4 | (30, 0.1) | (15, 3) | (100, 2.8) | 0.6 | 1.2 | 1.5 |

### *5.1.* Optimization via metaheuristic approaches

Metaheuristic is a search technique designed to find an appropriate response to an optimization problem that is difficult and time-consuming to solve. In the present study, DE and PSO metaheuristics are employed to maximize the non-linear profit function.

For three CE plans, both DE and PSO methods are used to compute the optimal values of the decision variables $(T_0, \xi_1, \xi_2, G, W_r)$ and the accompanying optimal profit. To handle constraints in our highly non-linear optimization model, we use penalty method for constraint handling is used along with metaheuristics (Lagaros et al., 2023).

### *5.1.1.* Differential evolution (DE)

DE consists three major processes mutation, crossover, and selection. For *mutation*, we initially select three different vectors $x_a, x_b$ and $x_c$ randomly and then create a donor vector $v_i^t$ using the following different mutation strategies after taking a differential weight $F \in [0,2]$ as a parameter and $R \in (0,1)$ as a random number (Jain & Singh, 2024). We shall use two mutation schemes viz., DE-1 and DE-2 in our current model to implement DE as follows:

DE/rand-to-best/1 or DE-1: $\quad v_i^{t+1} = x_i^t + R(x_{best}^t - x_i^t) + F(x_a^t - x_b^t)$

DE/Current-to-rand/1 or DE-2: $\quad v_i^{t+1} = x_i^t + R(x_a^t - x_i^t) + F(x_b^t - x_c^t)$

*Crossovers* are conducted via the binomial technique on each of the decision variables controlled by crossover probability $P_c = 0.8$.

In DE, the trial vectors produced from the crossover operation are selected as members of the current generation based on their fitness. To get the target population of individuals (best suited), the three fundamental activities mutation, crossover and selection, are carried out in each generation. We have run numerical experiments to choose suitable values of crossover probability and scale factor within an appropriate range. The input parameter values for DE are as follows:

**DE:** Scale factor $(F)$=0.6, $PopSize$=50, $MaxItr$=100.

### *5.1.2.* Particle swarm optimization (PSO)

In PSO, each particle follows a piecewise route that may be represented as a positional vector. Consider $N$ particles; the location of particle $i$ at iteration $t$ is $X_i^t$. For each particle, in addition to the position, we also have a velocity $V_i^t$. The location and velocity of each particle are modified in the next iteration as (Jain et al. 2022)

$$X_i^{t+1} = X_i^t + V_i^{t+1} \tag{38}$$
$$V_i^{t+1} = m_o V_i^t + c_1 r_1 (Pbest_i - X_i^t) + c_2 r_2 (Gbest - X_i^t) \tag{39}$$

Here, $r_1, r_2 \in (0,1)$; $m_o$ is inertia weight; $c_1, c_2$ are acceleration coefficients; $Pbest$ is the best position found by any particle $i$ and $Gbest$ is the best position among all the particles in swarm. To implement PSO, we set $c_1 = 2$, $c_2 = 2$, $PopSize$=50, $MaxItr$=300 as input parameters in our model.

Table 3 summarizes the maximum profit, average, and standard deviation for DE and PSO metaheuristics for the three carbon policies. The comparison of metaheuristic techniques DE-1, DE-2 and PSO based on statistical indices has also done to select the most suitable metaheuristic of our profit maximization problem. The convergence graphs of the suggested metaheuristic techniques are shown in Fig. 3.

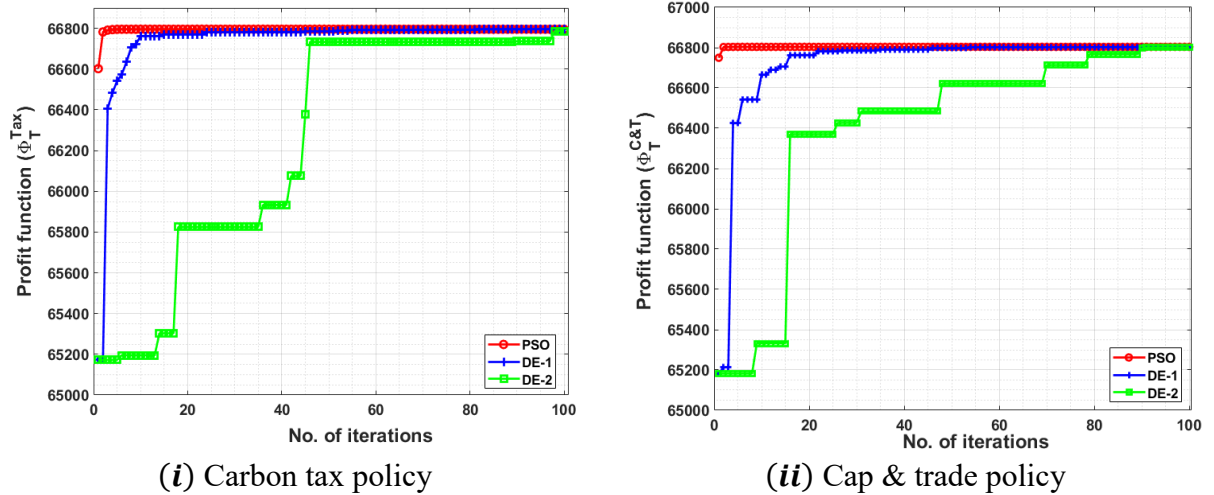

($i$) Carbon tax policy      ($ii$) Cap & trade policy

**Fig. 3:** Convergence graphs for DE-1, DE-2 and PSO in case of ($i$) Carbon tax policy, ($ii$) Cap & trade policy

The observations related to findings by the metaheuristics are as follows:

- In Table 4, we optimize the fitness function in case of three CE reduction policies. By tabulated values, we can decide which of the three algorithms are the best.
- The optimal profit, average value, and standard deviation of each technique are obtained. It is noticed that the PSO outperforms in comparison to DE-1 and DE-2 approaches.
- The convergence graphs are shown in Fig. 3 for the three metaheuristic techniques. It seems that the PSO converges more quickly than the DE-1 and DE-2 approaches.

**Table 3:** Comparative analysis of metaheuristics based on SD

| Carbon emission policies → | | Carbon Tax | Cap & Trade | Limited Emission |
|---|---|---|---|---|
| **DE-1** (DE/Rand-to-best/1) | **Max joint profit** | 66795.16 | 66803.57 | 23820.49 |
| | **Mean ($\bar{X}$)** | 66418.76 | 66800.96 | 22380.35 |
| | **SD ($\sigma$)** | 300.60 | 4.40 | 642.74 |
| **DE-2** (DE/Current-to-rand /1) | **Max joint profit** | 66786.45 | 66801.17 | 22048.12 |
| | **Mean ($\bar{X}$)** | 66042.29 | 66764.79 | 10464.98 |
| | **SD ($\sigma$)** | 448.86 | 24.48 | 21064.30 |
| **PSO** | **Max joint profit** | 66795.32 | 66803.71 | 29052.55 |
| | **Mean ($\bar{X}$)** | 66795.32 | 66803.71 | 28996.31 |
| | **SD ($\sigma$)** | $3.30 \times 10^{-10}$ | $6.13 \times 10^{-9}$ | 0.56 |

**Table 4:** Performance of metaheuristics for optimal policies

| Policies | Algorithms | Optimal Decisions | | | | | Profit value |
|---|---|---|---|---|---|---|---|
| | | $T_0^*$ | $\xi_1^*$ | $\xi_2^*$ | $W_r^*$ | $G^*$ | ($\varphi_T^*$) |
| Carbon Taxation | DE-1 | 0.6623 | 167.8949 | 95.7032 | 292.26 | 8.7962 | 66795.16 |
| | DE-2 | 0.6511 | 183.8457 | 93.5531 | 292.03 | 2.5875 | 66786.45 |
| | PSO | 0.6626 | 167.8648 | 93.6745 | 292.28 | 7.7568 | 66795.32 |
| Cap & Trade | DE-1 | 0.6621 | 167.2385 | 91.6493 | 292.27 | 7.1705 | 66803.57 |
| | DE-2 | 0.6644 | 163.5539 | 97.9678 | 292.13 | 2.4581 | 66801.17 |
| | PSO | 0.6624 | 167.8509 | 93.6717 | 292.27 | 7.7617 | 66803.71 |
| Limited Emission | DE-1 | 0.0054 | 0 | 0.0654 | 166.9600 | 0.01 | 23820.49 |
| | DE-2 | 0.0054 | 6.2052 | 2 | 110 | 0.01 | 22048.12 |
| | PSO | 0.0138 | 48.8437 | 47.6889 | 256.48 | 1.9693 | 29052.55 |

### 5.2. Profit function

The profit function of carbon tax policy model is concave w.r.t. the decision parameters for the carbon taxing strategy, as shown in Figs. 4(i-v). The profit function is non-linear with regards to the decision parameters, viz., $T_0, \xi_1, \xi_2, G$ and $W_r$. At a time, the two parameters are changed simultaneously to show the profit function's nature via surface plot. It is seen that metaheuristic optimization technique provides maximum profit as supported by the concave characteristics of profit function.

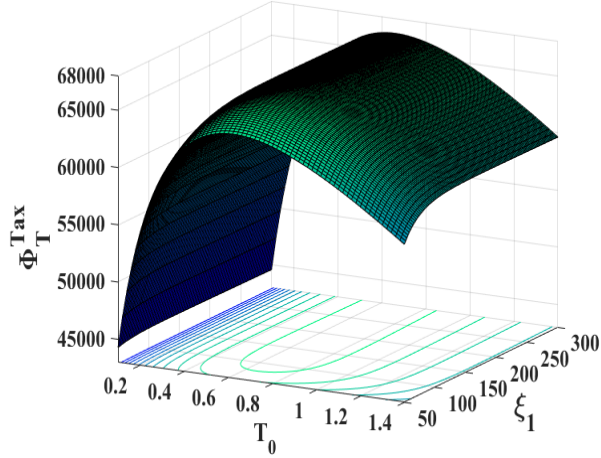
(i) Variation of $\varphi_T^{Tax}$ with $T_0$ and $\xi_1$

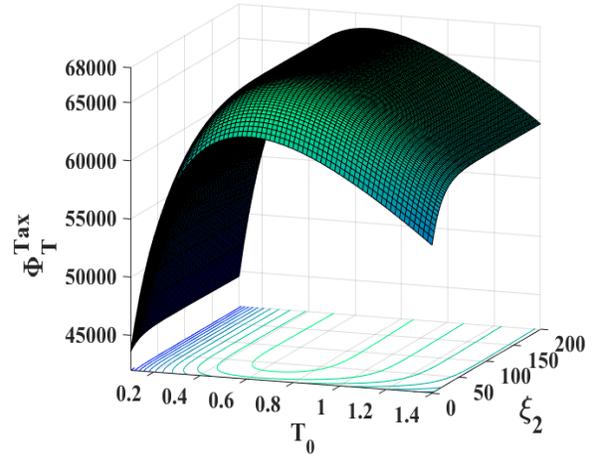
(ii) Variation of $\varphi_T^{Tax}$ with $T_0$ and $\xi_2$

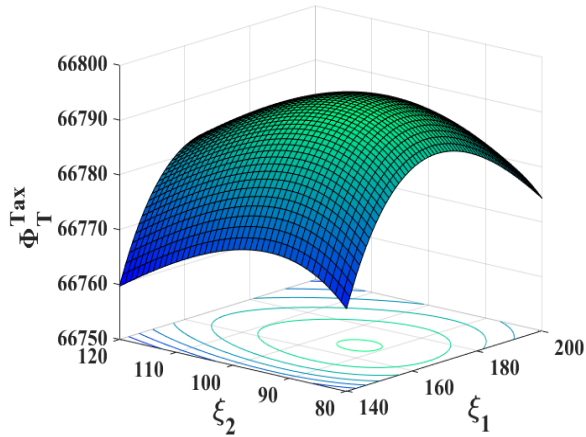
(iii) Variation of $\varphi_T^{Tax}$ with $\xi_1$ and $\xi_2$

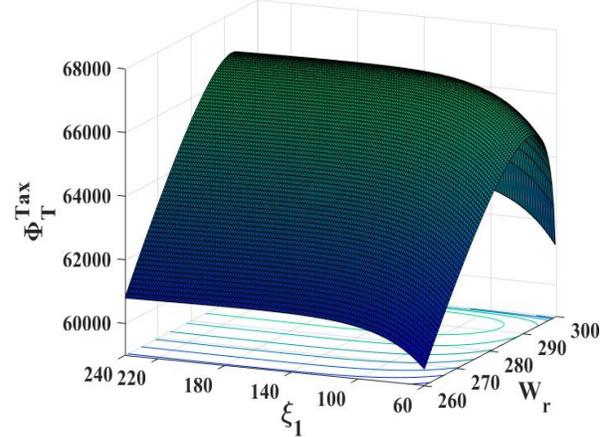
(iv) Variation of $\varphi_T^{Tax}$ with $\xi_1$ and $W_r$

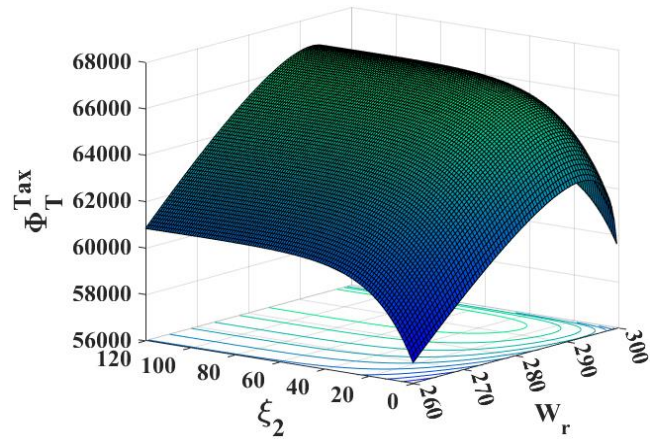
(v) Variation of $\varphi_T^{Tax}$ with $\xi_2$ and $W_r$

**Fig. 4:** Variation of joint profit function w.r.t. carbon taxation policy with (i) $T_0$ and $\xi_1$ (ii) $T_0$ and $\xi_2$ (iii) $\xi_1$ and $\xi_2$ (iv) $\xi_1$ and $W_r$ (v) $\xi_2$ and $W_r$

*5.3. ANFIS results*

To compare with findings from the metaheuristic PSO, the computational results using the ANFIS approach have been achieved. Using the neuro-fuzzy tool in MATLAB (R2023b) software, the ANFIS technique is now used to compute the results. The least square and backpropagation gradient descent methods are used to train the ANFIS training data set. The architecture of proposed model, including the number of nodes in each layer of the neural network, is elaborated in Table 5.

By adopting the linguistic values as very low, low, medium, high and very high of the respective input parameters, the membership function of each variable is shown in Fig. 5(i-ii). Fig. 6(i-vi) shows that ANFIS is given the approximately similar results by using training w.r.t. time and selling price. The demand function is price and stock dependent in our model.

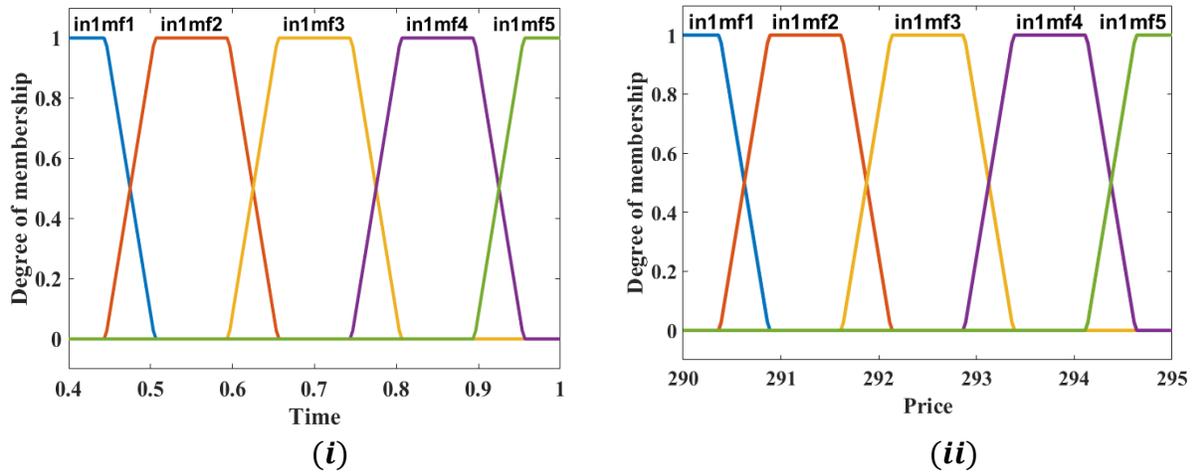

(*i*)            (*ii*)

**Fig. 5:** Membership functions of input variables (*i*) Time ($T_0$) and (*ii*) Price ($W_r$)

**Table 5:** ANFIS characteristics and Input parameters

| ANFIS Characteristics | Input Parameters | |
|---|---|---|
| | $T_0$ | $W_r$ |
| Type | Sugeno | Sugeno |
| No. of nodes | 24 | 24 |
| Number of linear and nonlinear parameters | 10, 20 | 10, 20 |
| Number of training data pairs | 61 | 101 |
| Number of fuzzy rules | 5 | 5 |
| AND method | Prod | Prod |
| OR method | Max | Max |
| Agg method | Sum | Sum |
| Defuzz method | wtaver | wtaver |

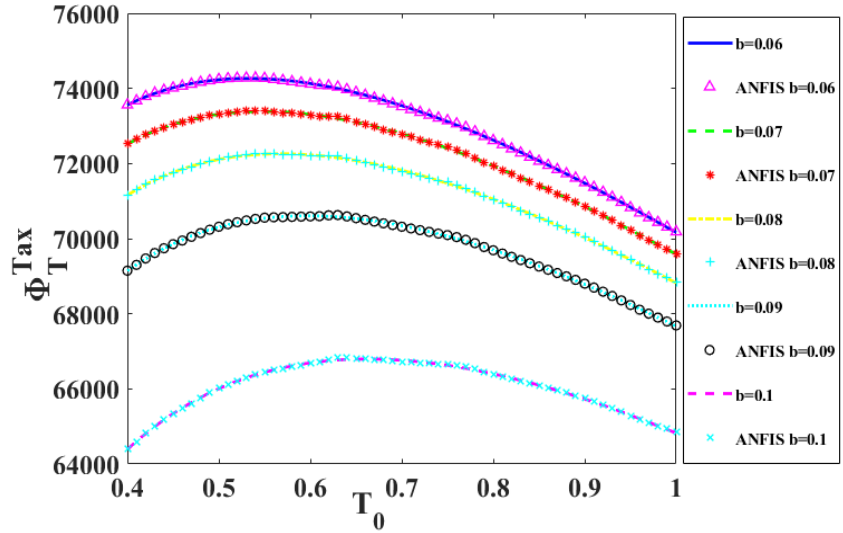

**Fig. 6(i):** Effect on $\varphi_T^{Tax}$ of $T_0$ by varying $b$ value

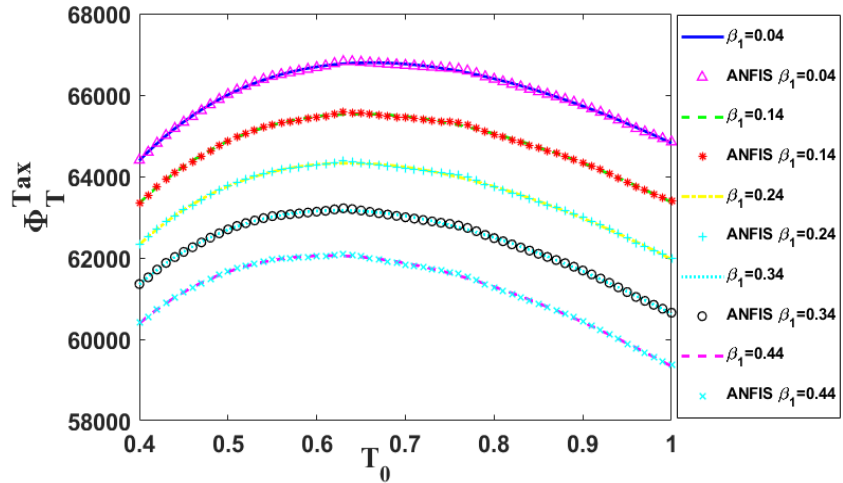

**Fig. 6(ii):** Effect on $\varphi_T^{Tax}$ of $T_0$ by varying $\beta_1$ value

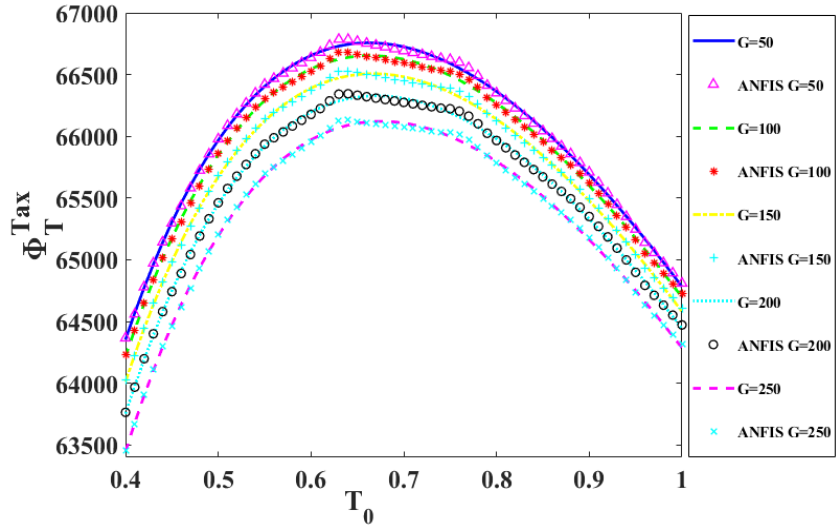

**Fig. 6(iii):** Effect on $\varphi_T^{Tax}$ of $T_0$ by varying $G$ value
24

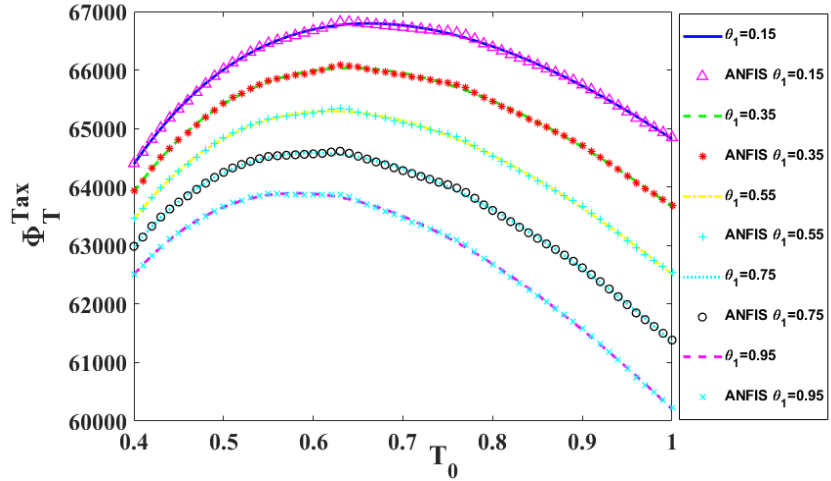

**Fig. 6(iv):** Effect on $\varphi_T^{Tax}$ of $T_0$ by varying $\theta_1$ value

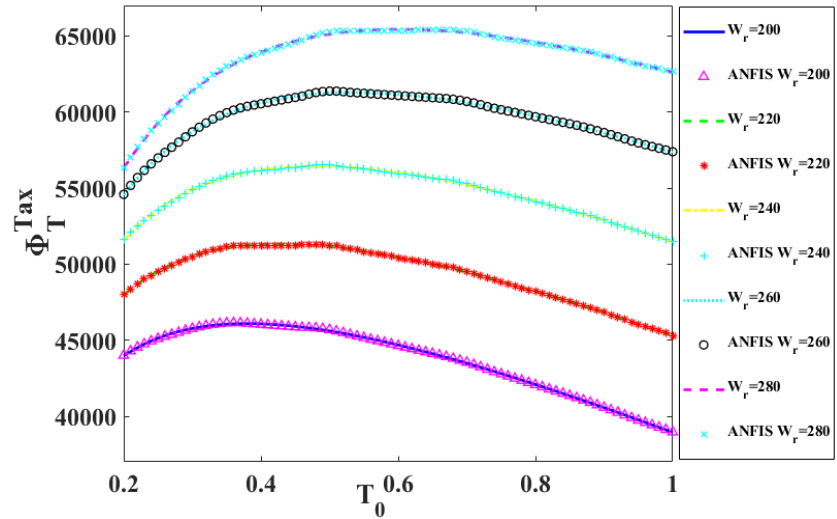

**Fig. 6(v):** Effect on $\varphi_T^{Tax}$ of $T_0$ by varying $W_r$ value

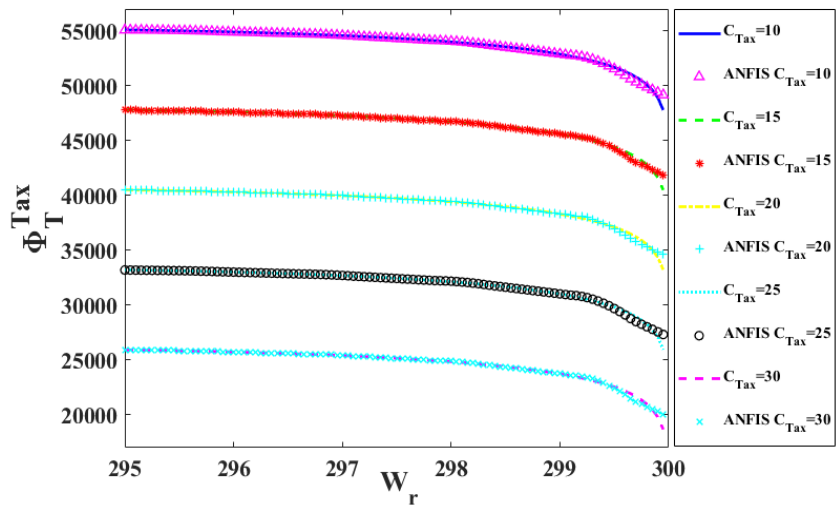

**Fig. 6(vi):** Effect on $\varphi_T^{Tax}$ of $W_r$ by varying $C_{Tax}$ value

## 6. Sensitivity analysis

We are interested to investigate the behaviour of the suggested model by examining the impacts of changes in specific parameters on the outcomes of the numerical simulations. The carbon emission parameters, cost parameters, production, depreciation, rework, holding parameters and many others are being assessed. The numerical results summarized in Tables 6-8 for carbon taxation policy only to demonstrate the sensitivity facilitate the following observations. There are three carbon emission strategies, each influenced by various parameters that directly impact a firm's joint profit. Here, we focus on examining how these parameters and their characteristics affect overall profit, individual profits of supply chain members, and green investments aimed at reducing emissions by considering carbon taxation policy.

### (I) Table 6

- A decrement in selling price ($W_m$) of the manufacturer negatively affects the joint profit function. An increment in $W_m$ shows enhancement in overall profit as well as the manufacturer profit.
- The joint profit ($\varphi_T^{Tax}$) seems too sensitive to the manufacturing cost ($C_p$) and rework cost ($C_r$). We must maintain the impact of manufacturing and rework costs on the profits of each stakeholder in a production SC.
- The joint profit ($\varphi_T^{Tax}$) decreases as we increase the emission costs w.r.t. production ($E_p$), transport ($E_t$) and holding ($E_{h1}$) of items. Some products spread pollution during production, transportation and storage before sale. Businessmen invest some money to sustain in the market with competitors to reduce this emission. Our model suggests that if we raise this investment, our earnings will be affected.
- The manufacturer invests in holding costs $h_p$ and $h_d$ w.r.t. perfect and deteriorating items, respectively. The holding cost of retailer ($h_r$) is less sensitive to the joint profit ($\varphi_T^{Tax}$). We see that $h_p$ is more sensitive to joint profit ($\varphi_T^{Tax}$) than $h_d$ and $h_r$ w.r.t. the percent changes in both costs. When there is 40% decrement in $h_p$, the joint profit shows 8.66% increment. On the other side, 40% increment in $h_p$ shows 6.69% decrement in the joint profit. Thus, we conclude that more investments to hold items can affect the profit for each player of SC.

### (II) Table 7

- Table 7 shows that the increment in production rate ($P$) and rework rate ($P_r$) of the manufacturer's cycle exhibits positive impact on the joint profit ($\varphi_T^{Tax}$). The production rate ($P$) is highly sensitive for time span $T_0$. It indicates that whenever a production system will increase the rework process of defectives and production process for the items, the joint profit as well as the profit of the producer increases.

**Table 6:** Sensitivity analysis of cost parameters

| Parameters | % Change | Decision variables | | | | | | | $\varphi_T^{Tax*}$ ($) | % change in $\varphi_T^{Tax*}$ |
|---|---|---|---|---|---|---|---|---|---|---|
| | | $T_0^*$ | $\xi_1^*$ | $\xi_2^*$ | $W_r^*$ | $G^*$ | $Z_m^*$ | $Z_r^*$ | | |
| $W_m$ | -40 | 0.7212 | 170.7275 | 93.7843 | 291.67 | 6.6707 | -6719.13 | 70791.42 | 64072.28 | -4.0767 |
| | -20 | 0.6925 | 169.3495 | 93.7374 | 291.97 | 7.1921 | -124.89 | 65534.11 | 65409.23 | -2.0751 |
| | 0 | 0.6626 | 167.8651 | 93.6741 | 292.28 | 7.7565 | 6493.11 | 60302.21 | 66795.32 | 0.0000 |
| | 20 | 0.6315 | 166.2558 | 93.5882 | 292.60 | 8.3703 | 13137.71 | 55098.57 | 68236.28 | 2.1573 |
| | 40 | 0.5989 | 164.5001 | 93.4744 | 292.92 | 9.0389 | 19812.32 | 49926.69 | 69739.01 | 4.4070 |
| $C_p$ | -40 | 0.6627 | 165.1448 | 93.6740 | 292.28 | 7.7574 | 8746.66 | 60300.82 | 69047.47 | 3.3717 |
| | -20 | 0.6627 | 166.5424 | 93.6741 | 292.28 | 7.7572 | 7619.80 | 60301.56 | 67921.36 | 1.6858 |
| | 0 | 0.6626 | 167.8651 | 93.6741 | 292.28 | 7.7565 | 6493.11 | 60302.21 | 66795.32 | 0.0000 |
| | 20 | 0.6626 | 169.1219 | 93.6742 | 292.28 | 7.7559 | 5366.55 | 60302.79 | 65669.34 | -1.6857 |
| | 40 | 0.6626 | 170.3177 | 93.6745 | 292.28 | 7.7565 | 4240.09 | 60303.33 | 64543.43 | -3.3713 |
| $C_r$ | -40 | 0.6626 | 167.7962 | 93.6746 | 292.28 | 7.7568 | 6556.94 | 60302.18 | 66859.12 | 0.0955 |
| | -20 | 0.6626 | 167.8310 | 93.6738 | 292.28 | 7.7568 | 6525.02 | 60302.2 | 66827.22 | 0.0478 |
| | 0 | 0.6626 | 167.8651 | 93.6741 | 292.28 | 7.7565 | 6493.11 | 60302.21 | 66795.32 | 0.0000 |
| | 20 | 0.6626 | 167.8984 | 93.6740 | 292.28 | 7.7566 | 6461.19 | 60302.22 | 66763.42 | -0.0478 |
| | 40 | 0.6626 | 167.9330 | 93.6738 | 292.28 | 7.7561 | 6429.29 | 60302.22 | 66731.52 | -0.0955 |
| $E_p$ | -40 | 0.6626 | 167.8157 | 93.6741 | 292.28 | 7.7567 | 6540.44 | 60302.17 | 66842.61 | 0.0708 |
| | -20 | 0.6626 | 167.8399 | 93.6741 | 292.28 | 7.7562 | 6516.78 | 60302.18 | 66818.96 | 0.0354 |
| | 0 | 0.6626 | 167.8651 | 93.6741 | 292.28 | 7.7565 | 6493.11 | 60302.21 | 66795.32 | 0.0000 |
| | 20 | 0.6626 | 167.8899 | 93.6738 | 292.28 | 7.7565 | 6469.46 | 60302.21 | 66771.67 | -0.0354 |
| | 40 | 0.6626 | 167.9143 | 93.6746 | 292.28 | 7.7567 | 6445.80 | 60302.23 | 66748.03 | -0.0708 |
| $E_t$ | -40 | 0.6626 | 167.9299 | 93.6741 | 292.28 | 7.7567 | 7359.05 | 60302.24 | 67661.29 | 1.2965 |
| | -20 | 0.6626 | 167.8974 | 93.6740 | 292.28 | 7.7568 | 6926.08 | 60302.22 | 67228.30 | 0.6482 |
| | 0 | 0.6626 | 167.8651 | 93.6741 | 292.28 | 7.7565 | 6493.11 | 60302.21 | 66795.32 | 0.0000 |
| | 20 | 0.6626 | 167.8323 | 93.6738 | 292.28 | 7.7564 | 6060.15 | 60302.19 | 66362.33 | -0.6482 |
| | 40 | 0.6626 | 167.7988 | 93.6742 | 292.28 | 7.7569 | 5627.17 | 60302.17 | 65929.35 | -1.2965 |
| $E_{h1}$ | -40 | 0.6711 | 168.1422 | 93.7514 | 292.35 | 7.5881 | 6571.75 | 60477.4 | 67049.15 | 0.3800 |
| | -20 | 0.6669 | 168.0032 | 93.7129 | 292.32 | 7.6728 | 6532.26 | 60389.57 | 66921.83 | 0.1894 |
| | 0 | 0.6626 | 167.8651 | 93.6741 | 292.28 | 7.7565 | 6493.11 | 60302.21 | 66795.32 | 0.0000 |
| | 20 | 0.6585 | 167.7280 | 93.6355 | 292.24 | 7.8389 | 6454.33 | 60215.28 | 66669.60 | -0.1882 |
| | 40 | 0.6544 | 167.5921 | 93.5975 | 292.20 | 7.9214 | 6415.87 | 60128.81 | 66544.68 | -0.3752 |
| $h_p$ | -40 | 0.8980 | 174.4056 | 95.3534 | 293.93 | 3.8456 | 8389.31 | 64195.48 | 72584.79 | 8.6675 |
| | -20 | 0.7599 | 170.8246 | 94.4728 | 293.06 | 5.9591 | 7343.69 | 62127.88 | 69471.57 | 4.0066 |
| | 0 | 0.6626 | 167.8651 | 93.6741 | 292.28 | 7.7565 | 6493.11 | 60302.21 | 66795.32 | 0.0000 |
| | 20 | 0.5896 | 165.3158 | 92.9338 | 291.56 | 9.2786 | 5783.59 | 58653.42 | 64437.02 | -3.5306 |
| | 40 | 0.5322 | 163.0608 | 92.2386 | 290.89 | 10.580 | 5181.86 | 57141.13 | 62322.99 | -6.6956 |
| $h_d$ | -40 | 0.6632 | 167.8581 | 93.6789 | 292.28 | 7.7462 | 6498.05 | 60313.19 | 66811.24 | 0.0238 |
| | -20 | 0.6629 | 167.8616 | 93.6764 | 292.28 | 7.7511 | 6495.57 | 60307.7 | 66803.28 | 0.0119 |
| | 0 | 0.6626 | 167.8651 | 93.6741 | 292.28 | 7.7565 | 6493.11 | 60302.21 | 66795.32 | 0.0000 |
| | 20 | 0.6624 | 167.8683 | 93.6720 | 292.28 | 7.7622 | 6490.64 | 60296.72 | 66787.36 | -0.0119 |
| | 40 | 0.6621 | 167.8717 | 93.6693 | 292.27 | 7.7666 | 6488.18 | 60291.22 | 66779.41 | -0.0238 |
| $h_r$ | -40 | 0.6636 | 167.9306 | 93.7120 | 292.27 | 7.7372 | 6472.15 | 60473.26 | 66945.41 | 0.2247 |
| | -20 | 0.6631 | 167.8977 | 93.6931 | 292.27 | 7.7464 | 6482.64 | 60387.72 | 66870.36 | 0.1123 |
| | 0 | 0.6626 | 167.8651 | 93.6741 | 292.28 | 7.7565 | 6493.11 | 60302.21 | 66795.32 | 0.0000 |
| | 20 | 0.6621 | 167.8321 | 93.6549 | 292.28 | 7.7662 | 6503.60 | 60216.69 | 66720.29 | -0.1123 |
| | 40 | 0.6616 | 167.7985 | 93.6352 | 292.29 | 7.7759 | 6514.09 | 60131.19 | 66645.28 | -0.2246 |

**Table 7:** Sensitivity analysis for the parameters $P, P_r, \theta_1, \theta_2, v_1, v_2, \beta_1, \beta_2$ and $\omega$

| Parameters | % Change | Decision variables | | | | | | | $\varphi_T^{Tax*}$ ($) | % change in $\varphi_T^{Tax*}$ |
| --- | --- | --- | --- | --- | --- | --- | --- | --- | --- | --- |
| | | $T_0^*$ | $\xi_1^*$ | $\xi_2^*$ | $W_r^*$ | $G^*$ | $Z_m^*$ | $Z_r^*$ | | |
| $P$ | -40 | 1.0695 | 169.1867 | 93.5473 | 292.26 | 7.8823 | 6284.91 | 58370.78 | 64655.69 | -3.2033 |
| | -20 | 0.8183 | 168.3954 | 93.6265 | 292.27 | 7.8045 | 6413.47 | 59563.09 | 65976.57 | -1.2258 |
| | 0 | 0.6626 | 167.8651 | 93.6741 | 292.28 | 7.7565 | 6493.11 | 60302.21 | 66795.32 | 0.0000 |
| | 20 | 0.5567 | 167.4851 | 93.7061 | 292.28 | 7.7238 | 6547.30 | 60805.24 | 67352.54 | 0.8342 |
| | 40 | 0.4800 | 167.2012 | 93.7287 | 292.28 | 7.7008 | 6586.55 | 61169.73 | 67756.28 | 1.4387 |
| $P_r$ | -40 | 0.6534 | 168.3458 | 93.6066 | 292.24 | 7.8646 | 6403.53 | 59719.33 | 66122.86 | -1.0067 |
| | -20 | 0.6591 | 168.0501 | 93.6486 | 292.26 | 7.7975 | 6459.24 | 60082.28 | 66541.52 | -0.3800 |
| | 0 | 0.6626 | 167.8651 | 93.6741 | 292.28 | 7.7565 | 6493.11 | 60302.21 | 66795.32 | 0.0000 |
| | 20 | 0.6650 | 167.7383 | 93.6911 | 292.29 | 7.7290 | 6515.88 | 60449.73 | 66965.62 | 0.2550 |
| | 40 | 0.6667 | 167.6459 | 93.7034 | 292.29 | 7.7097 | 6532.24 | 60555.57 | 67087.80 | 0.4379 |
| $\theta_1$ | -40 | 0.6701 | 155.3369 | 93.7426 | 292.34 | 7.6087 | 6575.22 | 60456.71 | 67031.93 | 0.3542 |
| | -20 | 0.6664 | 162.4070 | 93.7087 | 292.31 | 7.6822 | 6533.23 | 60379.28 | 66912.50 | 0.1754 |
| | 0 | 0.6626 | 167.8651 | 93.6741 | 292.28 | 7.7565 | 6493.11 | 60302.21 | 66795.32 | 0.0000 |
| | 20 | 0.6589 | 172.3036 | 93.6403 | 292.24 | 7.8296 | 6454.30 | 60225.48 | 66679.77 | -0.1730 |
| | 40 | 0.6553 | 176.0374 | 93.6062 | 292.21 | 7.9023 | 6416.46 | 60149.09 | 66565.55 | -0.3440 |
| $\theta_2$ | -40 | 0.6626 | 167.8665 | 85.1608 | 292.28 | 7.7562 | 6492.87 | 60312.62 | 66805.50 | 0.0152 |
| | -20 | 0.6626 | 167.8651 | 89.9553 | 292.28 | 7.7565 | 6493.00 | 60306.9 | 66799.90 | 0.0069 |
| | 0 | 0.6626 | 167.8651 | 93.6741 | 292.28 | 7.7565 | 6493.11 | 60302.21 | 66795.32 | 0.0000 |
| | 20 | 0.6626 | 167.8641 | 96.7120 | 292.28 | 7.7564 | 6493.25 | 60298.13 | 66791.38 | -0.0059 |
| | 40 | 0.6626 | 167.8639 | 99.2811 | 292.28 | 7.7570 | 6493.38 | 60294.49 | 66787.88 | -0.0111 |
| $v_1$ | -40 | 0.6623 | 258.4528 | 93.6663 | 292.27 | 7.7727 | 6405.41 | 60282.62 | 66688.03 | -0.1606 |
| | -20 | 0.6625 | 202.8472 | 93.6713 | 292.28 | 7.7622 | 6459.22 | 60294.86 | 66754.08 | -0.0617 |
| | 0 | 0.6626 | 167.8651 | 93.6741 | 292.28 | 7.7565 | 6493.11 | 60302.21 | 66795.32 | 0.0000 |
| | 20 | 0.6627 | 143.6904 | 93.6760 | 292.28 | 7.7529 | 6516.58 | 60307.09 | 66823.66 | 0.0424 |
| | 40 | 0.6628 | 125.9193 | 93.6778 | 292.28 | 7.7489 | 6533.83 | 60310.59 | 66844.41 | 0.0735 |
| $v_2$ | -40 | 0.6625 | 167.8555 | 141.9258 | 292.28 | 7.7589 | 6495.36 | 60244.03 | 66739.39 | -0.0837 |
| | -20 | 0.6626 | 167.8612 | 112.4413 | 292.28 | 7.7577 | 6493.95 | 60279.76 | 66773.71 | -0.0323 |
| | 0 | 0.6626 | 167.8651 | 93.6741 | 292.28 | 7.7565 | 6493.11 | 60302.21 | 66795.32 | 0.0000 |
| | 20 | 0.6627 | 167.8669 | 80.5951 | 292.28 | 7.7556 | 6492.56 | 60317.7 | 66810.26 | 0.0224 |
| | 40 | 0.6627 | 167.8689 | 70.9171 | 292.28 | 7.7559 | 6492.16 | 60329.08 | 66821.24 | 0.0388 |
| $\beta_1$ | -40 | 0.6650 | 167.7003 | 93.6904 | 292.29 | 7.7298 | 6544.73 | 60457.35 | 67002.08 | 0.3096 |
| | -20 | 0.6638 | 167.7829 | 93.6828 | 292.28 | 7.7432 | 6518.87 | 60379.65 | 66898.52 | 0.1545 |
| | 0 | 0.6626 | 167.8651 | 93.6741 | 292.28 | 7.7565 | 6493.11 | 60302.21 | 66795.32 | 0.0000 |
| | 20 | 0.6615 | 167.9461 | 93.6658 | 292.27 | 7.7694 | 6467.44 | 60225.03 | 66692.47 | -0.1540 |
| | 40 | 0.6603 | 168.0262 | 93.6578 | 292.27 | 7.7824 | 6441.88 | 60148.1 | 66589.98 | -0.3074 |
| $\beta_2$ | -40 | 0.6624 | 167.8757 | 93.6735 | 292.28 | 7.7587 | 6490.87 | 60289.99 | 66780.86 | -0.0217 |
| | -20 | 0.6625 | 167.8707 | 93.6730 | 292.28 | 7.7578 | 6491.99 | 60296.1 | 66788.09 | -0.0108 |
| | 0 | 0.6626 | 167.8651 | 93.6741 | 292.28 | 7.7565 | 6493.11 | 60302.21 | 66795.32 | 0.0000 |
| | 20 | 0.6627 | 167.8595 | 93.6749 | 292.28 | 7.7556 | 6494.25 | 60308.31 | 66802.55 | 0.0108 |
| | 40 | 0.6628 | 167.8539 | 93.6752 | 292.28 | 7.7545 | 6495.38 | 60314.41 | 66809.79 | 0.0217 |
| $\omega$ | -40 | 0.6627 | 167.8672 | 93.6745 | 292.28 | 7.2739 | 6491.91 | 60302.96 | 66794.87 | -0.0007 |
| | -20 | 0.6626 | 167.8650 | 93.6748 | 292.28 | 7.8457 | 6492.88 | 60302.32 | 66795.20 | -0.0002 |
| | 0 | 0.6626 | 167.8651 | 93.6741 | 292.28 | 7.7565 | 6493.11 | 60302.21 | 66795.32 | 0.0000 |
| | 20 | 0.6626 | 167.8653 | 93.6741 | 292.28 | 7.1260 | 6492.81 | 60302.41 | 66795.22 | -0.0001 |
| | 40 | 0.6627 | 167.8679 | 93.6746 | 292.28 | 6.1651 | 6492.17 | 60302.80 | 66794.97 | -0.0005 |

**Table 8:** Sensitivity analysis for the parameters $l_1$, $l_2$, $f_d$, $f_r$, $b$, $\eta$, $C_{Tax}$, $d_1$ and $i_c$

| Parameters | % Change | Decision variables | | | | | $Z_m^*$ | $Z_r^*$ | $\varphi_T^{Tax*}$ ($) | % change in $\varphi_T^{Tax*}$ |
|---|---|---|---|---|---|---|---|---|---|---|
| | | $T_0^*$ | $\xi_1^*$ | $\xi_2^*$ | $W_r^*$ | $G^*$ | | | | |
| $l_1$ | -40 | 0.6629 | 167.8817 | 93.6773 | 292.28 | 0.4157 | 6484.80 | 60307.73 | 66792.52 | -0.0042 |
| | -20 | 0.6628 | 167.8772 | 93.6760 | 292.28 | 2.8168 | 6487.15 | 60305.98 | 66793.14 | -0.0033 |
| | 0 | 0.6626 | 167.8650 | 93.6741 | 292.28 | 7.7565 | 6493.11 | 60302.21 | 66795.32 | 0.0000 |
| | 20 | 0.6622 | 167.8441 | 93.6703 | 292.27 | 15.5193 | 6504.27 | 60295.96 | 66800.23 | 0.0074 |
| | 40 | 0.6616 | 167.8118 | 93.6646 | 292.27 | 26.3159 | 6522.32 | 60286.82 | 66809.14 | 0.0207 |
| $l_2$ | -40 | 0.6620 | 167.8281 | 93.6673 | 292.27 | 24.1781 | 6511.43 | 60289.89 | 66801.32 | 0.0090 |
| | -20 | 0.6624 | 167.8538 | 93.6721 | 292.28 | 12.7425 | 6498.68 | 60298.46 | 66797.14 | 0.0027 |
| | 0 | 0.6626 | 167.8650 | 93.6741 | 292.28 | 7.7565 | 6493.11 | 60302.21 | 66795.32 | 0.0000 |
| | 20 | 0.6627 | 167.8703 | 93.6753 | 292.28 | 5.1714 | 6490.23 | 60304.14 | 66794.37 | -0.0014 |
| | 40 | 0.6628 | 167.8740 | 93.6755 | 292.28 | 3.6709 | 6488.56 | 60305.27 | 66793.82 | -0.0022 |
| $f_d$ | -40 | 0.6654 | 167.6672 | 93.6938 | 292.29 | 7.7254 | 6556.09 | 60486.05 | 67042.14 | 0.3695 |
| | -20 | 0.6640 | 167.7660 | 93.6842 | 292.28 | 7.7413 | 6524.54 | 60393.94 | 66918.47 | 0.1844 |
| | 0 | 0.6626 | 167.8650 | 93.6741 | 292.28 | 7.7565 | 6493.11 | 60302.21 | 66795.32 | 0.0000 |
| | 20 | 0.6612 | 167.9619 | 93.6636 | 292.27 | 7.7720 | 6461.84 | 60210.83 | 66672.67 | -0.1836 |
| | 40 | 0.6599 | 168.0571 | 93.6544 | 292.27 | 7.7871 | 6430.70 | 60119.81 | 66550.52 | -0.3665 |
| $f_r$ | -40 | 0.6645 | 167.9830 | 93.6910 | 292.29 | 7.7195 | 6454.70 | 60340.54 | 67039.04 | 0.3649 |
| | -20 | 0.6636 | 167.9241 | 93.6830 | 292.29 | 7.7386 | 6473.90 | 60321.4 | 66917.16 | 0.1824 |
| | 0 | 0.6626 | 167.8650 | 93.6741 | 292.28 | 7.7565 | 6493.11 | 60302.21 | 66795.32 | 0.0000 |
| | 20 | 0.6617 | 167.8050 | 93.6657 | 292.27 | 7.7747 | 6512.35 | 60282.95 | 66673.51 | -0.1824 |
| | 40 | 0.6608 | 167.7459 | 93.6568 | 292.26 | 7.7930 | 6531.59 | 60263.65 | 66551.75 | -0.3646 |
| $b$ | -40 | 0.3817 | 145.7886 | 81.4348 | 175.03 | 14.2108 | 12334.68 | 22315.68 | 34650.35 | -48.1246 |
| | -20 | 0.5379 | 159.0912 | 88.6178 | 233.54 | 10.3683 | 9086.36 | 40967.74 | 50054.10 | -25.0635 |
| | 0 | 0.6626 | 167.8650 | 93.6741 | 292.28 | 7.7565 | 6493.11 | 60302.21 | 66795.32 | 0.0000 |
| | 20 | 0.7689 | 174.4828 | 97.6113 | 351.16 | 5.8490 | 4283.97 | 80044.81 | 84328.78 | 26.2495 |
| | 40 | 0.8627 | 179.8255 | 100.840 | 410.15 | 4.4056 | 2334.06 | 100066.4 | 102400.44 | 53.3048 |
| $\eta$ | -40 | 0.7632 | 172.4389 | 99.5929 | 289.89 | 5.3079 | 4403.37 | 54540.16 | 58943.53 | -11.7550 |
| | -20 | 0.7076 | 169.9909 | 96.3607 | 291.31 | 6.6341 | 5559.70 | 57951.79 | 63511.48 | -4.9163 |
| | 0 | 0.6626 | 167.8650 | 93.6741 | 292.28 | 7.7565 | 6493.11 | 60302.21 | 66795.32 | 0.0000 |
| | 20 | 0.6257 | 166.0113 | 91.3819 | 292.99 | 8.7164 | 7259.85 | 62049.7 | 69309.55 | 3.7641 |
| | 40 | 0.5948 | 164.3778 | 89.3840 | 293.54 | 9.5491 | 7902.63 | 63415.76 | 71318.39 | 6.7715 |
| $C_{Tax}$ | -40 | 0.6730 | 168.2278 | 93.7733 | 292.37 | 1.1713 | 7487.82 | 60537.83 | 68025.65 | 1.8419 |
| | -20 | 0.6678 | 168.0460 | 93.7241 | 292.32 | 4.5325 | 6989.72 | 60419.75 | 67409.47 | 0.9195 |
| | 0 | 0.6626 | 167.8650 | 93.6741 | 292.28 | 7.7565 | 6493.11 | 60302.21 | 66795.32 | 0.0000 |
| | 20 | 0.6575 | 167.6846 | 93.6248 | 292.23 | 10.4401 | 5997.42 | 60185.5 | 66182.92 | -0.9168 |
| | 40 | 0.6525 | 167.5062 | 93.5754 | 292.19 | 12.6257 | 5502.42 | 60069.66 | 65572.08 | -1.8313 |
| $d_1$ | -40 | 0.6626 | 167.9297 | 93.6743 | 292.28 | 7.7564 | 7359.05 | 60302.24 | 67661.29 | 1.2965 |
| | -20 | 0.6626 | 167.8975 | 93.6744 | 292.28 | 7.7564 | 6926.09 | 60302.21 | 67228.30 | 0.6482 |
| | 0 | 0.6626 | 167.8650 | 93.6741 | 292.28 | 7.7565 | 6493.11 | 60302.21 | 66795.32 | 0.0000 |
| | 20 | 0.6626 | 167.8327 | 93.6738 | 292.28 | 7.7566 | 6060.15 | 60302.19 | 66362.33 | -0.6482 |
| | 40 | 0.6626 | 167.7991 | 93.6740 | 292.28 | 7.7565 | 5627.19 | 60302.16 | 65929.35 | -1.2965 |
| $i_c$ | -40 | 0.6626 | 167.1677 | 93.6742 | 292.28 | 7.7567 | 7094.01 | 60301.85 | 67395.86 | 0.8991 |
| | -20 | 0.6626 | 167.5184 | 93.6741 | 292.28 | 7.7566 | 6793.55 | 60302.04 | 67095.59 | 0.4495 |
| | 0 | 0.6626 | 167.8650 | 93.6741 | 292.28 | 7.7565 | 6493.11 | 60302.21 | 66795.32 | 0.0000 |
| | 20 | 0.6626 | 168.2060 | 93.6740 | 292.28 | 7.7561 | 6192.70 | 60302.35 | 66495.05 | -0.4495 |
| | 40 | 0.6626 | 168.5428 | 93.6737 | 292.28 | 7.7563 | 5892.26 | 60302.53 | 66194.79 | -0.8991 |

- It is realized that the joint profit ($\varphi_T^{Tax}$) as well as individual profit of the manufacturer ($Z_m$) and retailer ($Z_r$) increase when the deterioration rates for the manufacturer ($\theta_1$) and retailer ($\theta_2$) decrease.
- Our model suggests that if any business can control over the amount of the degraded items, there may a boost in overall profit of SC.
- The preservation constants for the manufacturer ($v_1$) and retailer ($v_2$) are positively correlated with the joint profit ($\varphi_T^{Tax}$) but do not significantly affect $\varphi_T^{Tax}$. An increment in $v_1$ and $v_2$ leads to higher quantity of items and less deteriorated goods.
- In any business, loss will increase as it becomes more likely that non-faulty items discarded as defective. We can conclude that a 40% increment in type-I error ($\beta_1$) gives a fall up to 0.31% in the joint profit ($\varphi_T^{Tax}$).
- If a faulty item is declared as non-defective at the manufacturer cycle, increment in the probability that a faulty item can sell as non-defective ($\beta_2$), shows a very little increment in individual profit ($Z_m$) as well as overall profit ($\varphi_T^{Tax}$).
- Individual profit of each stakeholder as well as the collective profit ($\varphi_T^{Tax}$) have negative impact of a proportion ($\omega$) of green investments.

### (III) Table 8

- If the retailer will receive goods with IE of type-II ($\beta_2$), he needs to invest in goodwill impairment with a fraction ($f_r$) otherwise he will suffer goodwill loss as a result of selling defective products. We can see that the increment in $f_r$ negatively affects the profit functions of each chain member ($Z_m$ and $Z_r$) as well joint profit ($\varphi_T^{Tax}$).
- Stock dependent demand parameter 'η' is positively correlated with the joint profit ($\varphi_T^{Tax}$). When there is 40% increment and decrement in η, $\varphi_T^{Tax}$ shows 6.77% increment as well as 11.75% decrement.
- The price dependent demand parameter ($b$) is extremely sensitive to both the joint profit ($\varphi_T^{Tax}$) and the profitability of the chain members individually ($Z_m$ and $Z_r$). It indicates that increment in $b$ will affect the demand of the customers and price of the items. The profit of the manufacturer ($Z_m$) lowers down by the increment in $b$. The profit of the retailer ($Z_r$) shows increment w.r.t. '$b$'. The joint profit ($\varphi_T^{Tax}$) shows 53.30% rise when '$b$' goes up to 40% increment. The $\varphi_T^{Tax}$ decreases up to 48.12% when there is fall of 40% in '$b$'.
- In any centralized green SC, every chain member wants to invest less costs and earns more profit. If the emission of harmful gases increases by any business firm, government applies the more amount of carbon taxes ($C_{Tax}$). It can affect the overall profit. According to our model, if $C_{Tax}$ rises by 40%, overall profit ($\varphi_T^{Tax}$) decreases to 1.83%.
- If the manufacturer needs to deliver orders at distance '$d_1$' and if this distance will increase up to 40%, the profit of the manufacturer ($Z_m$) decreases. So, there is a decrement of 1.29% in the joint profit.

- We can see that $f_d$ and $i_c$ do not significantly affect the joint profit ($\varphi_T^{Tax}$) of the SC. When these parameters rise by 40%, overall profit ($\varphi_T^{Tax}$) decreases to 0.36% and 0.89%, respectively. Also, when $l_1$ rises by 40%, overall profit ($\varphi_T^{Tax}$) increases to 0.02%.

## 7. Managerial insights

Strategic approaches based on decision-making are required to stay ahead in the competitive world of management. Managerial insights are the foundation of effective leadership which give executives and managers a perspective through which firms can grab opportunities, and advance the performance of their organizations. The following are some key insights explained:

a) To optimize profitability in a carbon-regulated environment, managers should concentrate on cutting production costs, streamlining holding and rework procedures, and funding emission-reduction technology.
b) When choosing the best carbon emission strategy, managers should take the industry's legislation and specific circumstances into consideration.
c) The joint profit is positively impacted by increasing production and rework rates. To optimize profitability, managers should make investments in improving the effectiveness of rework procedures and production efficiency.
d) Profit is lowered by higher emission costs. It's imperative to invest in more environment friendly processes and technologies that reduce emissions during production, shipping, and storage.
e) The total profit increases dramatically when the selling price rises. Selling prices should be carefully chosen by managers to guarantee that they are high enough to optimize profit without adversely affecting demand.
f) For some parameters, the ANFIS methodology yields results similar to those of metaheuristic techniques, indicating that it may be a valuable tool for demand prediction and price optimization in uncertain environments. Under dynamic market situations, managers can use ANFIS to make data-driven decisions.

## 8. Conclusions

With available green energy options, reducing CE is an important strategy which can be used to mitigate climate change's effects and maintain a green environment. Following three distinct CE regulatory policies, viz., carbon taxes, cap & trade, and restricted emission, our study focused on the centralised SC with manufacturing process and emission reduction impacts by using green investment. This article incorporated the realistic features of deterioration, preservation, imperfect production, rework, IE, price and stock-dependent demand, CE and green investment. Two metaheuristics DE with different mutation schemes and PSO have employed to solve the non-linear constrained optimization problems of three CE policies.

When making future decisions and improving two-tier SC inventory models with inspection errors, it may be helpful to consider how sensitive the different performance measures that have been established in the current model. The proposed green SC inventory model is validated using numerical simulation and optimization. The sensitivity analysis done can be used to examine the parameter variation impact on the profit function. The outcomes of the green SC study demonstrate the benefits of cap & trade and carbon taxing systems over restricted emission policy. It is noticed that the demand function parameters that are stock and price dependent are quite sensitive. In order to decrease shortages and waste while increasing SC profit, the firm should focus on pricing and stock management. To increase profitability, businesses must first cut manufacturing, rework, storage, and emissions expenses. To reduce deterioration, the firm/industries are advised to invest in preservation technology. The developed model investigates how preservation efforts affect the rate of spoiling and motivated decision-makers to opt the best preservation investment. The rework can help to eliminate waste caused by imperfect production. Based on our studied model concludes that the manufacturer's selling price should be controlled so that it cannot negatively impact the profit function.

In order to verify the viability of utilizing fuzzy parameters and neural networks in our sustainable supply chain inventory model, ANFIS technique is deployed effectively. The validation of ANFIS in our model provides insight into its potential application in a number of intricate supply chain inventory situations where closed-form analytical results are difficult to derive. By predicting demand and lowering carbon emissions, the sensitivity of various factors used may also be beneficial to increase the profit. These techniques also offer useful insights for forecasting the demand, quantity, shortages, and profit. Our study can be further extended by considering stochastic demand, fuzzy membership for imprecise parameters, trade credit policy, etc.

**Data Availability Statement**

No additional data have been used in this study.

**Acknowledgement**

The author would like to thank the All India Council for Technical Education (AICTE).

**Declaration of Interest**

The authors declare that they have no known competing financial interests or personal relationships that could have appeared to influence the work reported in this paper.